\documentclass[11pt]{amsart}
\usepackage{amssymb,amscd,verbatim,mathrsfs}
\usepackage{a4wide}
\usepackage[matrix,arrow]{xy}

\usepackage{color}

\usepackage[dvips]{graphicx}
\usepackage{epsfig,enumerate}
\numberwithin{equation}{section}
\newtheorem{theorem}{Theorem}[section]
\newtheorem{proposition}[theorem]{Proposition}
\newtheorem{lemma}[theorem]{Lemma}

\theoremstyle{definition}
\newtheorem{definition}[theorem]{Definition}
\theoremstyle{remark}

\DeclareMathOperator{\Tr}{Tr}
\newcommand{\Dd}{{D}}

\newcommand{\Vv}{\mathcal{V}}

\newcommand{\la}{\lambda}

\newcommand{\Hh}{\mathcal{H}}

\newcommand{\pT}{{}^\Phi\!T}
\newcommand{\pVv}{{}^\Phi\!\Vv}
\newcommand{\cz}{\mathbb C}
\newcommand{\rz}{\mathbb R}
\newcommand{\zz}{\mathbb Z}
\newcommand{\Cc}{\mathcal{C}}

\newcommand{\Nn}{\mathcal{N}}

\def\tr{\text{tr}}

\def\sR{\hbox{I\kern-.1667em\hbox{R}}}
\newcommand{\R}{\mathbb R}

\newcommand{\Z}{\mathbb Z}

\def\D{{\mathcal D}}

\def\tr{\hbox{Tr}}

\def\sl{\mathrm{SL}(2,\mathbb{R})}

\def\cusLchi{L^2(\Gamma\backslash G, \chi)_{\mathrm{cus}}}

\def\cusLchi{L^2(\Gamma\backslash G,\chi)_{\mathrm{cus}}}

\def\ctLchi{L^2(\Gamma\backslash G, \chi)_{\mathrm{ct}}}
\def\resLchi{L^2(\Gamma\backslash G,\chi)_{\mathrm{res}}}

\def\sLchi{L^2(\Gamma\backslash G, \chi)}
\def\disLchi{L^2(\Gamma\backslash G, \chi)_{\mathrm{dis}}}

\def\t{\mathrm{t}}
\def\tr{\mathrm{tr}}
\def\cst{\mathrm{cst}}
\def\ct{\mathrm{ct}}
\def\dis{\mathrm{dis}}
\def\Tr{\mathrm{Tr}}
\def\pr{\mathrm{pr}}
\def\p{\mathrm{p}}
\def\d{\mathrm{d}}
\begin{document}

\title[Adiabatic limit of the Eta invariant]
{Adiabatic limit of the Eta invariant over cofinite quotient of
$\mathrm{PSL}(2,\R)$}

\author{Paul Loya}
\address{Department of Mathematics \\Binghamton University\\
Binghamton\\NY 13902\\ U.S.A. } \email{paul@math.binghamton.edu}

\author{Sergiu Moroianu}
\address{Institutul de Matematic\u a al Academiei Rom\^ane\\ P.O. Box
1-764\\ RO-014700 Bucha\-rest\\ Romania}
\email{moroianu@alum.mit.edu}

\author{Jinsung Park}
\address{School of Mathematics\\ Korea Institute for Advanced Study\\
207-43\\ Cheongnyangni 2-dong\\ Dong\-daemun-gu\\ Seoul 130-722\\
Korea } \email{jinsung@kias.re.kr}

\thanks{2000 Mathematics Subject Classification.
Primary: 58J28, 58J50. Secondary: 11F72, 22E46.}

\date{\today}

\begin{abstract}
{We study the adiabatic limit of the eta invariant of the Dirac
operator over cofinite quotient of $\mathrm{PSL}(2,\R)$, which is
a \emph{noncompact} manifold with a \emph{nonexact} fibred-cusp
metric near the ends.}
\end{abstract}

\maketitle

%%%%%%%%%%%%%%%%%%%%%%%%%%%%%%%%%%%%%%%%%%%%%%%%%%%%%%%%%%%%%%%%%%%%
\section{Introduction}
%%%%%%%%%%%%%%%%%%%%%%%%%%%%%%%%%%%%%%%%%%%%%%%%%%%%%%%%%%%%%%%%%%%%%

The eta invariant of a Dirac operator was introduced in the
seminal paper \cite{APS75} as the boundary correction term in the
Atiyah-Patodi-Singer index formula. Since the paper of Atiyah
\emph{et al.}\ \cite{APS75}, the eta invariant has found a place
in many branches of mathematics. One particular aspect of the eta
invariant which has found fruitful applications is the study of
the \emph{adiabatic limit} of the eta invariant, in which the eta
invariant on the total space of a fibration is investigated when
the fiber is collapsed. This was initiated by Witten \cite{W} and
later proved independently by Bismut and Freed \cite{BF} and
Cheeger \cite{Ch87}. Expanding on the earlier work of
\cite{BF,Ch87,W}, Bismut and Cheeger \cite{Bis-ch89} and then Dai
\cite{Dai91} studied the adiabatic limit for a general fibration
of compact manifolds. We also refer to the recent work of Moroianu
\cite{Mor04}, who analyzed the adiabatic limit for general
families of first order elliptic operators from the view point of
the calculus of adiabatic pseudodifferential operators.

In this paper we study the spectral properties of the Dirac
operator and compute the adiabatic limit of its eta invariant on a
certain three-dimensional noncompact manifold $X$, which is given
by a cofinite quotient of $\mathrm{PSL}(2,\R)$. Then $X$ is a
circle bundle that fibers over a \emph{non-compact} Riemann
surface $\Sigma$ with cusps,
\begin{equation} \label{S}
\xymatrix{{S}^1  \ar[r] &  X \ar[d] \\
              & \Sigma_{g,\kappa} \, , }
\end{equation}
where $g,\kappa$ denote the genus and the number of cusps,
respectively, of the base Riemann surface which has finite volume.
There has been much interest in more general spectral problems for
the case when the fibers are circles over compact bases
\cite{Amm-Bar98,Be00,Dai-Zhang95,Hitchin,NiL98,NicL99,Seade-Steer,Z}.
However, in all these papers the base manifold is compact. In this
paper we consider the spectral properties and the eta invariant in
the case when the base is not compact.

Replacing the circle ${S}^1$ in \eqref{S} with a circle of radius
$r$, then choosing a spin structure, we denote the corresponding
Dirac operator by $\Dd_r$. The first purpose of this paper is to
study the spectral properties of $D_r$ defined by non-exact fibred
cusp metrics near the ends. In particular, we analyze the
dependence of the continuous spectrum of $D_r$ on the choice of
spin structure. Postponing terminology concerning trivial
spin structures to Section 2, we can now state our first result:

\begin{theorem}\label{t:cont-spec} With $\kappa^{\t}$ denoting
the number of ends with trivial spin structures, the continuous
spectrum of the Dirac operator $\Dd_r$ consists of
$\kappa^\t$-copies of countably many families given by
\[
\left(\, -\infty, \, -\frac{r}{2}-|m|
\big(1+r^{-2}\big)^{1/2}\,\right]\,\bigcup\, \left[\,-\frac{r}{2}+
|m|\big(1+r^{-2}\big)^{1/2},\, +\infty\,\right)
\]
indexed by odd integers $m\in 1+2\mathbb{Z}$ if the spin structure
along the $S^1$-fiber is trivial, or by even integers
$m\in2\mathbb{Z}$ otherwise.
\end{theorem}

This theorem can be regarded as a generalization of the result of
B\"ar in \cite{Bar00} to the fibred cusp case where the continuous
spectrum of $D_r$ depends on the spin structures of the
$S^1$-fibers and of the $S^1$ cross sections of the base manifold
$\Sigma_{g,\kappa}$ near the ends. Another novelty of this theorem
is that the Riemannian metric over the cusps are \emph{not}
{exact} fibred cusp metrics, which have been extensively studied
in, for example, \cite{Lau-mor05,LMP,Vaill}. It is because of the
non-exact fibred cusp metrics that the continuous spectrum of
$\Dd_r$ is quite complicated.

The second main result of this paper is the adiabatic limit of the
eta invariant of $\Dd_r$ as the fiber is collapsed (that is, $r
\to 0$). As we observed in Theorem \ref{t:cont-spec}, the Dirac
operator $\Dd_r$ has continuous spectrum as well as discrete
spectrum; moreover, the corresponding odd heat kernel of $\Dd_r^2$
is not trace class. Therefore, the standard definition of the eta
function using the eigenvalues or the trace of the odd heat kernel
are not valid in our situation. This requires us to define a
``regularized'' eta invariant, which is reminiscent of the $b$-eta
invariant of Melrose \cite{Mel93} and similar to a regularized eta
invariant used by Park \cite{Park05} to analyze eta invariants on
hyperbolic manifolds with cusps. With $\eta(\Dd_r, s)$ denoting
the eta function of $\Dd_r$ defined through a regularized trace
similar to Melrose's $b$-trace \cite{Mel93} (see Definition in
\eqref{e:reg-tr}), the following is our main result:

\begin{theorem}\label{t:main thm} We assume that the spin
structure along the $S^1$-fiber is trivial.
\begin{enumerate}
\item For sufficiently small $r>0$, $\eta(\Dd_r,s)$ defined for
$\Re(s)\gg 0$ has a meromorphic extension over $\mathbb{C}$ and
may have a double pole at $s=1$ and simple poles on $
-\mathbb{N}\cup \{0,2\}$. Moreover, for a totally nontrivial spin
structure, $\eta(\Dd_r,s)$ may have only simple poles at
$-\mathbb{N}\cup\{0,1,2\}$.

\item For the eta invariant,
\[
\eta(\Dd_r):= \mathrm{Reg}_{s = 0} \, \eta(\Dd_r,s)
\]
where $\mathrm{Reg}_{s = 0}$ means to take the regular value at $s
= 0$, the following equality holds:
\begin{equation}\label{e:main}
\lim_{r\to 0}\eta(\Dd_r)=
-\frac{1}{12\pi}\mathrm{Vol}(\Sigma_{g,\kappa})=\frac16
\left(2-2g-{\kappa}\right)
\end{equation}
where $\mathrm{Vol}(\Sigma_{g,\kappa})$ is given w.r.t. the
Poincar\'e metric.
\end{enumerate}
\end{theorem}

For the compact case, and for the trivial spin structure, a result
corresponding to the formula \eqref{e:main} in Theorem \ref{t:main
thm} was proved by Seade and Steer \cite{Seade-Steer}, who also
obtained the original value of the eta invariant by applying the
APS index formula for a manifold with smooth boundary. In our
noncompact case, obtaining the original value of the eta invariant
will require an index formula for manifolds whose boundaries are
manifolds with non-exact fibred cusp ends. This problem will be
considered elsewhere.

The paper is organized as follows.  In Section \ref{sec-Dirac} we
develop the required background material, including a discussion
of spin structures and the Dirac operator $\Dd_r$.  In Section
\ref{sec-Anafibredcusp} we analyze the Dirac operator in terms of
the fibred cusp calculus of Mazzeo--Melrose \cite{Mazz-Mel98} and
we prove Theorem \ref{t:cont-spec}. In Section
\ref{sec-regtracandgeomside} we define the regularized eta
invariant and in Sections \ref{sec-FT},
\ref{sec-etafnprincipalseriespart}, and
\ref{sec-etafnDiscreteseriespart} we analyze the geometric and
spectral sides of the Selberg trace formula in our context, which
will be used to prove Theorem \ref{t:main thm}.

%%%%%%%%%%%%%%%%%%%%%%%%%%%%%%%%%%%%%%%%%%%%%%%%%%%%%%%%%%%%%%%%%%%%%%%%%
\section{Dirac operator and Spin structure}
\label{sec-Dirac}
%%%%%%%%%%%%%%%%%%%%%%%%%%%%%%%%%%%%%%%%%%%%%%%%%%%%%%%%%%%%%%%%%%%%%%%%%

In this section, we define the Dirac operator over the cofinite
quotient of $\mathrm{PSL}(2,\R)$ by {a} discrete subgroup.
Equivalently we consider the Lie group $G=\mathrm{SL}(2,\R)$ and a
discrete subgroup $\Gamma \subset G$ containing $\{\pm 1\}$, so
that the quotient $\Gamma\backslash G$ is the same as the quotient
{$(\Gamma/\{\pm 1\})\backslash\mathrm{PSL}(2,\R)$.}

For $r\in (0,\infty)$ we define a family of metrics $g_r$ over
$G$ such that the left translations of $E:=r^{-1}C, A, H$
are orthonormal w.r.t. $g_r$ where  $C,A,H$ is a basis of
$\mathfrak{g}=\mathrm{sl}(2,\R)$ given by
\begin{equation}\label{def:basis}
C=\begin{pmatrix} 0 & 1\\ -1 & 0 \end{pmatrix},\quad
A=\begin{pmatrix} 0 &1 \\ 1 & 0 \end{pmatrix},\quad
H=\begin{pmatrix} 1 & 0 \\ 0& -1\end{pmatrix}.
\end{equation}
{Recall} that the corresponding Levi-Civita connection $\nabla^r$
is determined by {the Koszul formula}
\begin{multline*}
2 g_r(\nabla^r_X Y, Z) = X g_r(Y,Z) + Y g_r(Z,X) - Z
g_r(X,Y)\\
+g_r([X,Y],Z) -g_r([X,Z],Y) - g_r([Y,Z],X)
\end{multline*}
where $X,Y,Z$ denote vector fields over $G$.

 Since $G$ is topologically the same as ${S}^1\times
\mathcal{H}$ where $\mathcal{H}$ is the Poincar\'e upper half
plane, there are two spin structures on $G$. We choose the one
determined by the left invariant trivialization. Denoting the
lifted connection to spinor bundle by the same notation
$\nabla^r$, we define the Dirac operator by
\[
\widehat{D}_r:=E.\nabla_E^r + A.\nabla_A^r + H.\nabla_H^r
\]
where $X.$ denotes the Clifford action by $X$. By a
straightforward computation as in \cite{Hitchin},
\cite{Seade-Steer}, we obtain
\[
\widehat{D}_r \psi = \frac{1}{2} \big(\frac{2+r^2}{r^2}-2\big) C.
A. H. \psi
\]
for a basic spinor $\psi$.

We twist $\widehat{D}_r$ by multiplying the volume element
$\omega:=E.A.H.$ to define $\widetilde{D}_r$, that is,
\[
\widetilde{D}_r:=E.A.H.\widehat{D}_r,
\]
which has the following simplified form,
\[
\widetilde{D}_r\psi = \Big(\frac{2-r^2}{2r}\Big)\, \psi
\]
for a basic spinor $\psi$. The Clifford algebra generated by
$E,A,H$ has the \emph{Pauli matrix representation} given by
\begin{equation*}
E\mapsto\begin{pmatrix} i & 0 \\ 0 & -i \end{pmatrix}, \quad
A\mapsto\begin{pmatrix} 0 & 1 \\ -1 & 0 \end{pmatrix}, \quad
H\mapsto\begin{pmatrix} 0 & -i \\ -i & 0 \end{pmatrix}.
\end{equation*}
Then we have
\[
\omega E\mapsto -i \begin{pmatrix} 1 & 0 \\ 0 &-1\end{pmatrix},
\quad \omega A\mapsto \begin{pmatrix} 0 & 1 \\ -1 & 0
\end{pmatrix}, \quad \omega H \mapsto i\begin{pmatrix} 0 & 1\\ 1 & 0
\end{pmatrix}.
\]
From these, for any spinor $\alpha \psi_1+ \beta \psi_2$ with
basic spinors $\psi_1,\psi_2$ and smooth functions $\alpha,\beta$
on $G$, we have the following representation of $\widetilde{D}_r$,
\begin{equation}\label{rep}
\widetilde{D}_r \begin{pmatrix}\alpha\\ \beta\end{pmatrix} =
\Big(\frac{2-r^2}{2r} \Big)\,\begin{pmatrix}\alpha\\
\beta\end{pmatrix} +
\begin{pmatrix} -i E & A+i H \\ -A+i H & i E \end{pmatrix}
\begin{pmatrix}\alpha\\ \beta\end{pmatrix} .
\end{equation}
Now we let
\[
Z:=-iC, \qquad 2X_+:= A-iH,\qquad 2X_-:= A+iH
\]
(Our convention is slightly different from the one in
\cite{Seade-Steer}). Then they satisfy
\[
[Z,X_+]=2X_+, \qquad [Z,X_-]=-2X_-, \qquad [X_+,X_-]=Z,
\]
and we have
\[
\widetilde{D}_r = \Big(\frac{2-r^2}{2r}\Big)+
\begin{pmatrix} r^{-1}Z & 2X_-
\\ -2X_+ & -r^{-1}Z \end{pmatrix} \quad \text{acting on} \quad C^\infty(G)\oplus C^{\infty}(G) .
\]
It is also easy to check that
\begin{equation}\label{e:Lap}
\widetilde{D}_r^2= \begin{pmatrix} - (A^2+H^2+r^{-2}C^2) & 0 \\ 0
& - (A^2+H^2+r^{-2}C^2) \end{pmatrix} + \ \text{lower order
terms},
\end{equation}
hence the Dirac Laplacian $\widetilde{D}^2_r$ is a generalized
Laplacian whose principal symbol is given by the metric $g_r$.

To define the Dirac operator over $X=\Gamma\backslash G$, let us
discuss on the spin structures on $X=\Gamma\backslash G$, which
will play a crucial role throughout this paper. First recall that
there are $\vert H^1(X, \mathbb{Z}_2)\vert$-number of spin
structures over $X$ since every $3$-dimensional manifold is spin.
This can be understood from the following diagram,
\[ \begin{CD}
\widetilde{\mathcal{S}}  @>>>
\mathcal{S}\\
@VVV       @VVV \\
\widetilde{X} @> \pi >> X\\
\end{CD}
\]
where $\widetilde{\mathcal{S}}, \mathcal{S}$ are the Spin bundles
over the universal covering manifold $\widetilde{X}$ and $X$
respectively. Since $\widetilde{\mathcal{S}}\cong
\widetilde{X}\times \mathrm{Spin}(3)$ and
$\widetilde{\mathcal{S}}\cong\pi^*\mathcal{S}$, the possible Spin
bundle $\mathcal{S}$ is given by the $\mathbb{Z}_2$-representation
$\rho$ of $\pi_1(X)$ as follows:
\begin{equation}\label{e:spinstr}
\mathcal{S}_\rho = \widetilde{X}\times_{\rho} \mathrm{Spin}(3)
\end{equation}
with the obvious $\mathbb{Z}_2$-action to $\mathrm{Spin}(3)$.
Therefore, each $\mathbb{Z}_2$-representation of $\pi_1(X)$
provides us with inequivalent spin structure on $X$. Recall
\begin{multline*}
\pi_1(X)=\big\{\, x_i,y_i, h_j, k \ \ \big| \ \ 1\leq i\leq g,\,
1\leq
j \leq \kappa, \\
\prod^g_{i=1}[x_i,y_i]\prod^{\kappa}_{j=1} h_i=1, \
[x_i,k]=[y_i,k]=[h_j,k]=1 \, \big\},
\end{multline*}
where $g,\kappa$ denote the number of genus and cusps of the base
Riemann surface $\Sigma_{g,\kappa}$ of the fibration \eqref{S}.

Among spin structures, there are spin structures which are
determined by {those} $\mathbb{Z}_2$-{representations} $\rho$ of
$\pi_1(X)$ with $\rho(h_j)=-1$ for some $j$. Such a spin structure
over Riemann surface $\Sigma_{g,\kappa}$ is called
\emph{nontrivial along the cusp} corresponding to $j$ in
\cite{Bar00}. Following \cite{Bar00}, we call such a spin
structure \emph{nontrivial along the cusp} if $\rho(h_j)=-1$ for
the corresponding $j$, and \emph{totally nontrivial spin
structure} if it is determined by {a}
$\mathbb{Z}_2$-representation $\rho$ with
\begin{equation}\label{nontrivial}
{\text{$\rho(h_j)=-1$ for all $j=1,\ldots,\kappa$.}}
\end{equation}
From the relation of the generators of $\pi_1(X)$, there is the
following obstruction for this,
\[
\prod^\kappa_{j=1}\rho(h_j)=1.
\]
Hence, in this case, the number of cusps $\kappa$ should be even.
We {distinguish} two classes of spin structures according to
(non)triviality of spin structure along the fiber $S^1/\{\pm 1\}$.
We call the spin structure {\emph{trivial along the fiber} if the
spin structure is trivial along the fiber $S^1/\{\pm 1\}$ (or
equivalently, if the representation $\rho$ maps the generator $k$
to $1$), and \emph{nontrivial along the fiber} otherwise.} Note
that if the spin structure is trivial along the $S^1/\{\pm
1\}$-fiber, this spin structure does not extend to a spin
structure over the disc bundle over $\Sigma_{g,\kappa}$. From the
above discussion we have

\begin{proposition}\label{spinstr}There are
$2^{{2g+\kappa}}$ spin structures over $X=\Gamma\backslash G$.
There {exist} totally nontrivial spin structures over $X$ if and
only if $\kappa$ is even.
\end{proposition}

For the trivial representation of $\pi_1(X)$, the resulting Spin
bundle denoted by $\mathcal{S}_{1}$ is topologically trivial,
determined by the left invariant trivialization over
$X=\Gamma\backslash G$. The associated spinor bundle
$\Sigma_1=\mathcal{S}_{1}\times_{\mathrm{Spin}(3)} \Sigma(3)$
(where $\Sigma(3)\cong \mathbb{C}^2$ is the spinor representation
of $\mathrm{Spin}(3)\cong \mathrm{SU}(2)$) is therefore also
trivial and has the following relation with other $\Sigma_\rho$,
\[
\Sigma_\rho= \Sigma_1\otimes \mathbb{C}_\rho
\]
where $\mathbb{C}_\rho\to X$ is the flat line bundle associated to
$\rho$.

From the definition of $\widetilde{D}_r$ over $G$ and the equality
\eqref{e:spinstr}, we can see that the induced Dirac operator from
$\widetilde{D}_r$ pushed down to $X$ has the following form,
\begin{equation}\label{fd}
{D}_r = \Big(\frac{2-r^2}{2r}\Big)+
\begin{pmatrix} r^{-1}Z & 2X_-
\\ -2X_+ & -r^{-1}Z \end{pmatrix} \quad \text{acting on}
\quad C^\infty_0(\Gamma\backslash G, \chi)
\end{equation}
where $C^\infty_0(\Gamma\backslash G, \chi)$
($\chi=\rho\oplus\rho$) consists of the smooth functions with
co-compact supports such that $f(\gamma x)=\chi(\gamma)f(x)$ for
$\gamma\in\pi_1(X)$, $x\in G$. We also denote the $L^2$-completion
of $D_r$ (w.r.t. certain metric metric explained in
\eqref{e:Haar}) by $\Dd_r$, that is,
\begin{equation}\label{Dd}
\Dd_r \, : \ \sLchi \longrightarrow \sLchi.
\end{equation}

%%%%%%%%%%%%%%%%%%%%%%%%%%%%%%%%%%%%%%%%%%%%%%%%%%%%%%%%%%%%%%%%%%%%
\section{Analysis of fibred cusp operators}
\label{sec-Anafibredcusp}
%%%%%%%%%%%%%%%%%%%%%%%%%%%%%%%%%%%%%%%%%%%%%%%%%%%%%%%%%%%%%%%%%%%%%

In this section we show that the metrics $g_r$ are of
\emph{conformal fibred cusp} type. Consequently, we show
that the Dirac operators $\Dd_r$ belong to the class of weighted
fibred cusp operators introduced by Mazzeo and Melrose
\cite{Mazz-Mel98}, and we prove Theorem \ref{t:cont-spec}.

First we introduce some subgroups of $G=\sl$,
\begin{equation}\label{e:param}
N_0=\biggl\lbrace \begin{pmatrix} 1 & x\\ 0& 1\end{pmatrix}
\bigg\vert\ x\in\R\biggr\rbrace,    A_0=\biggl\lbrace
\begin{pmatrix} e^{\frac u2} & 0\\ 0&
e^{-\frac u2}\end{pmatrix} \bigg\vert\ u\in\R \biggr\rbrace, \\
K=\biggl\lbrace \begin{pmatrix} \cos\theta & \sin\theta\\
-\sin\theta &\cos\theta \end{pmatrix}\biggr\rbrace.
\end{equation}
Then the standard parabolic subgroup $P_0$ is given by $N_0A_0Z$
where $Z=\{\pm 1\}\subset K$ and any \emph{parabolic subgroup} $P$
of $G$ is conjugate to $P_0$ by an element {$k_P$} in $K$. A
parabolic subgroup $P$ has a decomposition $P=N_P A Z$ where $N_P$
is the derived group of $P$ and $A$ is any conjugate of $A_0$ in
$P$, to be called a \emph{Cartan subgroup}. It is clear that $A_0$
is the unique Cartan subgroup $P_0$ with Lie algebra orthogonal to
that of $K$. Therefore, $P$ has a unique Cartan subgroup with the
same property. From now on, we assume that the pair $(P,A)$
satisfies this condition. For such a pair $(P,A)$ with $N=N_P$, we
have the \emph{Iwasawa decomposition} $G=NAK$.

For a given $\Gamma\subset G$, a parabolic subgroup $P$ of $G$ is
called \emph{$\Gamma$-cuspidal} if $N=N_P$ contains a nontrivial
element of $\Gamma$. As one knows, the finitely many ends of
$X\cong \Gamma\backslash G$ are parametrized by $\Gamma$-conjugacy
classes $\{ P\}_\Gamma=\{ \gamma P\gamma^{-1} \, | \, \gamma \in
\Gamma/\Gamma_P\}$ where $\Gamma_P:=\Gamma\cap P$. Let $P$ be a
$\Gamma$-cuspidal parabolic subgroup of $G$ corresponding to one
end of $X=\Gamma\backslash G$. This subgroup determines a cusp
$c_P$, {an incomplete manifold} which is identified with a
neighborhood of {the cuspidal end} of the quotient
$\Gamma_P\backslash G$.

Assume first that $P=P_0$ is the standard parabolic subgroup of
$G$. The manifold  $\Gamma_P\backslash G$ has two commuting free
$S^1$ actions: the action of $K$ to the right and that of
$\Gamma_{N_0}\backslash N_0$ to the left where
$\Gamma_{N_0}:=\Gamma\cap N_0$. The first $S^1$ action exists in
fact globally on $X=\Gamma\backslash G$, while the second one
exists only on the cusp. Let $\gamma_l:=\begin{pmatrix}
1&l\\0&1\end{pmatrix}$ be the generator of $\Gamma_{N_0}$. Then
$\Gamma_P\backslash G$ is identified with $\rz/l\zz\times
\rz\times \rz/2\pi\rz$ by the map
\begin{equation}\label{coord}
(x,u,\theta)\mapsto \begin{pmatrix}
1&x\\0&1\end{pmatrix}\begin{pmatrix}
e^{\frac u2}&0\\0&e^{-\frac u2}\end{pmatrix}
\begin{pmatrix}  \cos\theta & \sin\theta\\
-\sin\theta &\cos\theta \end{pmatrix}.\end{equation} By projection
on the last two terms, we view this as the total space of a
fibration with fiber $S^1$. Note that this fibration makes
sense only near the end and is not the fibration in \eqref{S}
where the roles of two ${S}^1$'s are reversed.

As seen above, the spinor bundle corresponding to the
representation $\rho$ is the spinor bundle for the trivial
representation, twisted by the flat line bundle $\mathbb{C}_\rho$
defined by $\rho$. The Dirac operator on $G$ has been computed in
\eqref{fd} with respect to the orthonormal vector fields
$r^{-1}C$, $A$, $H$ defined in \eqref{def:basis} and the
representation $\chi$. The same expression holds on the spinor
bundle on the cusp $\Gamma_P\backslash G$, where the vector fields
$r^{-1}C$, $A$, $H$ now act on $\Sigma_\rho$. There is no
ambiguity about the action of these vector fields since the
twisting bundle $\mathbb{C}_\rho$ is flat.

Introduce the function $\nu :=e^{-u}$ on the cusp and glue the
``boundary at infinity'' $\rz/l\zz\times \{\nu =0\}\times
\rz/2\pi\rz$ to the cusp, thus getting a manifold with boundary
$\overline{\Gamma_P\backslash G}$. The $S^1$-fibration structures
extend to the boundary. We will show that for each fixed $r$, the
metric $g_r$ on $X$ is conformal to a fibred cusp metric (with
respect to the fibration of the boundary induced from the
$\Gamma_{N_0}\backslash N_0$ action). In the coordinates $(x,\nu
=e^{-u},\theta)$ of $\Gamma_P\backslash G$, the coefficients of a
matrix $\begin{pmatrix} a&b\\c&d\end{pmatrix}$ are given by the
inverse of the map \eqref{coord}:
\begin{align*}
x=\frac{ac+bd}{c^2+d^2},&& \nu
=c^2+d^2,&&\theta=-\arctan\left(\frac{c}{d}\right).
\end{align*}
We compute then
\begin{align}\label{e:diff}
&{E=\ \ r^{-1}}\partial_\theta,& \notag\\
&A=-\cos 2\theta\partial_\theta+2\nu^{-1}\cos 2\theta \partial_x
-2{\nu} \sin 2\theta \partial_{\nu} ,&\\
&H=\ \ \sin 2\theta \partial_\theta -2{\nu} ^{-1} \sin 2\theta
\partial_x -2{\nu}  \cos 2\theta \partial_{\nu} .&\notag
\end{align}
These equalities also can be found at p.52 of
\cite{Gel-Gr-Sh} or p.115 of \cite{Lang}. It follows that in the
coordinates $(x,{\nu} ,\theta)$ the metric $g_r$ is given by
\[\frac{1}{4{\nu} ^2}d{\nu} ^2+\frac{{\nu} ^2}{4}dx^2+r^2(d\theta
+\frac{\nu}{{2}}dx)^2,\]
thus
\[g_\Phi:=\frac{4}{{\nu} ^2}g_r=\left(\frac{d{\nu} }{{\nu}^2}\right)^2
+r^2\left(\frac{2d\theta}{{\nu}}+dx\right)^2+dx^2.\] This is what
is called a fibred cusp metric, or a smooth metric on the fibred
cusp tangent bundle. To define this, consider the subalgebra
$\pVv$ of the Lie algebra $\Vv$ of smooth vector fields on the
manifold with boundary $\overline{\Gamma_P\backslash G}$, spanned
over $C^\infty(\overline{\Gamma_P\backslash G})$ by the vector
fields
\begin{align*}
V_{\nu}:={\nu}^2\partial_{\nu},&&V_\theta
:={\nu}\partial_\theta,&&V_x:=\partial_x.
\end{align*}
This sub-algebra is by definition a free
$C^\infty(\overline{\Gamma_P\backslash G})$-module so it is the
space of sections of a smooth vector bundle over
$C^\infty(\overline{\Gamma_P\backslash G})$; this vector bundle,
called $\pT \Gamma_P\backslash G$, comes equipped with a bundle
morphism to the usual tangent bundle
$T\overline{\Gamma_P\backslash G}$, induced from the inclusion of
the spaces of sections $\pVv\hookrightarrow \Vv$, which is an
isomorphism over $\Gamma_P\backslash G$.

Since $\pVv$ is a Lie algebra and the metric $g_\Phi$ defined
above is non-degenerate and smooth on fibred cusp vector fields,
it follows immediately from the Cartan formula that the
Levi-Civita connection on $\Gamma_P\backslash G$ with respect to
the metric $g_\Phi$ extends to the boundary in the sense that for
every $V_i,V_j,V_k\in\pVv$, we have
\[\langle \nabla_{V_i}V_j,V_k\rangle\in
C^\infty(\overline{\Gamma_P\backslash G}).\] The spinor bundle
$\Sigma_\rho$ extends over the boundary, such that the Clifford
multiplication by $V_i$ is a smooth map. Take now the orthonormal
frame
\begin{align*}
V_1:=V_{\nu},&&V_2:=V_x-V_\theta/2,&&V_3:=\frac{1}{2r}
V_\theta.\end{align*} {Its relationship to the
global frame $(E,A,H)$ is deduced from \eqref{e:diff}:}
{\begin{equation}\label{trf}\begin{split}
V_1=& -\frac{\nu}{2}(\sin 2\theta A+\cos 2\theta H)\\
V_2=&\ \ \ \frac{\nu}{2}(\cos 2\theta A-\sin 2\theta H)\\
V_3=&\ \ \ \, \frac{\nu}{2}E.
\end{split}\end{equation}}

Denote by $V$ a local
lift to the spinor bundle of the orthonormal frame $(V_1,V_3,V_2)$. Let
$\sigma:\overline{\Gamma_P\backslash G}\to\Sigma(3)$ be a smooth
map into the $3$-spinor representation space. It follows from the
local formula
\begin{equation}\label{dphi}\begin{split}
D_\Phi[V,\sigma]=&\sum_{i=1}^3c(V_i)\left([V,V_i(\sigma)]+\frac{1}{2}\sum_{j<k}
c(V_j)c(V_k)\langle \nabla_{V_i}V_j,V_k\rangle\right)\\
=&\left(c(V_1)({\nu}^2\partial_{\nu}{-}\frac{\nu}{2})
+c(V_2)(\partial_x-\frac{{\nu}\partial_\theta}{2})+c(V_3)\partial_\theta\frac{{\nu}}{2r}
{-} r\frac{\nu}{4}\right)[V,\sigma]
\end{split}\end{equation}
that the Dirac operator with respect to $g_\Phi$ (defined first on
compactly supported spinors over $\Gamma_P\backslash G$) extends
to smooth spinors up to the boundary. Such an operator, a
combination of fibred cusp vector fields and of smooth bundle
endomorphisms down to the boundary $\{{\nu}=0\}$, is called a
fibred cusp differential operator. Thus
\[D_\Phi\in \mathrm{Diff}^1_\Phi(\overline{\Gamma_P\backslash G},\Sigma_\rho).\]

The Dirac operator changes very nicely with respect to conformal
changes of the metric. We simply have
\[D_r=2{\nu}^{-2}\circ D_\Phi\circ {\nu}\]
so for $r>0$, the Dirac operator $D_r$ is a differential operator
in ${\nu}^{-1}\mathrm{Diff}^1_\Phi(\overline{\Gamma_P\backslash
G},\Sigma_\rho)$.

The \emph{normal operator} $\Nn(D_\Phi)(\theta,\xi,\tau)$ of $D_\Phi$
(see \cite{Mazz-Mel98}) is obtained by replacing formally
\begin{align*}
V_{\nu}\mapsto i\xi,&& V_\theta\mapsto i\tau
\end{align*}
and then restricting to ${\nu}=0$. The result is a family of
differential operators on the {fibers of the boundary fibration
(the $x$-circles) with coefficients in the spinor bundle}, with
parameters $\theta\in S^1$, $(\xi,\tau)\in \rz^2$:
\[\Nn(D_\Phi)(\theta,\xi,\tau)=
c(V_1)i\xi+c(V_2)
\left(\partial_x-\frac{i\tau}{2}\right)+c(V_3)\frac{i\tau}{2r}.\]

\begin{definition}
The operator $D_\Phi$ is called \emph{fully elliptic} if
$\Nn(D_\Phi)(\theta,\xi,\tau)$ is invertible for all
$(\theta,\xi,\tau)\in S^1\times\rz^2$.
\end{definition}

If $D_\Phi$ is fully elliptic, then {by the results of
\cite{Mazz-Mel98} it has a parametrix inside the calculus of
fibered cusp pseudodifferential operators $\Psi_\Phi^{-1}(X)$}
modulo compact operators.

{\bf Proof of Theorem \ref{t:cont-spec}}. According to the
decomposition principle {(see e.g.\ \cite[Proposition 1]{Bar00})},
the essential spectrum of $\Dd_r$ can be computed outside a
compact subset of $X$, thus it is a superposition of the essential
spectra of any self-adjoint extension of $D_r$ over each {cuspidal
end $c_P$ defined by $\nu_P<\epsilon_P$. We must make sure that
such an extension exists (the Dirac operator on a manifold with
boundary may not admit self-adjoint extensions, e.g. on $\rz^+$).
We may take for instance the Atiyah-Patodi-Singer boundary
condition at the torus boundary $\{\nu_P=\epsilon_P\}$. Special
care is needed for the nullspace of the Dirac operator along the
torus, we allow in the domain only harmonic spinors of the form
$(u,c(V_3)u)$ where $u$ is in the $i$-eigenspace of $c(V_1)$}.

Since any $\Gamma$-parabolic subgroup $P$ is conjugated by an
element in the maximal compact subgroup $K$ to the standard
parabolic subgroup $P_0$, we see that the cusp corresponding to
$P$ is isometric to the ``canonical'' cusp $P_0\backslash G$. Thus
we can assume that we work with the canonical parabolic subgroup
$P_0$.

We have seen above that $D_r$ belongs to
${\nu}^{-1}\mathrm{Diff}^1_\Phi(\overline{\Gamma_P\backslash
G},\Sigma_\rho)$ {near the cuspidal end.}
\begin{lemma}
The {fibered cusp} operator $D_\Phi$ is fully elliptic on the cusp
$c_P$ if and only if the spin structure is non-trivial along
$c_P$.
\end{lemma}
\begin{proof}
We have computed above the normal operator
$\Nn:=\Nn(D_\Phi)({\theta},\xi,\tau)$. Clearly, $\Nn$ is an
elliptic self-adjoint operator on the circle in the variable $x$.
Therefore $\Nn$ is invertible if and only if $\Nn^2$ is. Now by
the anti-commutation of the Clifford variables,
\[\Nn^2=\xi^2+\frac{\tau^2}{4r^2}+(i\partial_x-\frac{\tau}{2})^2.\]
This family of operators is independent of $\theta$; it is
strictly positive (hence invertible) for $(\xi,\tau)\neq
0\in\rz^2$. For $\xi=\tau=0$, $\Nn=-\partial_x^2$, so $\ker(\Nn)$
is made of those spinors which are constant in $x$ in the
trivialization $V$ of the spinor bundle. For fixed $\theta$, such
spinors exist globally on the $x$-circle if and only if the local
lift $V$ satisfies $V_{x=l}=V_{x=0}$. Now the frame
$(V_1,V_2,V_3)$ is obtained from $(E,A,H)$ by the transformation
\eqref{trf} which is constant in $x$; thus the lift $V$ exists
globally around the cusp if and only if the lift of $(E,A,H)$
exists globally around the cusp, which is by definition equivalent
to the triviality of the spin structure around the cusp $c_P$.
\end{proof}
If $D_\Phi$ is fully elliptic, it follows from the above lemma and
from the general theory of fibred cups operators that
$D_r=2{\nu}^{-2}D_\Phi {\nu}$ has a parametrix $Q\in
{\nu}\Psi_\Phi^{-1}$ over the cusp $c_P$, modulo compact
operators. But $Q$ itself is compact due to the decaying weight
${\nu}$; hence the self-adjoint operator $\D_r$ has pure-point
spectrum over the cusp.

Conversely, assume that the spin structure is trivial along the
cusp. The operator $D_r=2{\nu}^{-2}D_\Phi {\nu}$ computed in
\eqref{dphi} has constant coefficients in $x$, thus it preserves
the orthogonal decomposition into zero-modes and high-energy modes
\[L^2(\Gamma_P\backslash G{\cap\{\nu_P<\epsilon_P\}}, \Sigma_\rho)=:
\Hh_0\oplus \Hh'\] where $\Hh_0$ is the space of {spinors}
constant in $x$ {in the trivialization $V$ (we have seen above
that $V$ exists globally around the cusp if the spin structure is
trivial along $c_P$) and $\Hh'$ its orthogonal complement}. Over
$\Hh'$, by the same argument as above, there exists a compact
parametrix of $D_r$ inside the fibered cusp calculus. Thus the
essential spectrum {of $D_r$ over the cusp $c_P$ only arises from
the zero-modes, i.e., it is the essential spectrum of the
operator}
\[2{\nu}^{-1}\left(c(V_1)({\nu}\partial_{\nu}{-}\frac12)
-c(V_2)\frac{\partial_\theta}{2} +c(V_3)\frac{\partial_\theta}{2r}
{-}\frac{r}{4}\right){\nu}\] acting in $L^2([0,\epsilon)\times
S^1,\Sigma_\rho,d{\nu} d\theta)$ {with any boundary condition at
$\epsilon$ which makes it self-adjoint.} We conjugate this
operator through the Hilbert space isometry
\begin{align*}L^2\left(d{\nu}d\theta\right)\to L^2
\left(\frac{d{\nu}}{{\nu}}d\theta\right)
&& \phi\mapsto {\nu}^\frac12 \phi.
\end{align*}
We get the operator
\[A_r=2c(V_1) {\nu}\partial_{\nu}
+\left({\frac{c(V_3)}{r}}-c(V_2)\right)\partial_\theta {-}\frac
r2.\] This can be again decomposed according to the frequencies in
the $\theta$ variable. {Note that although the local lift $V$ may
not exist globally, the ambiguity is locally constant so the
operator $i\partial_\theta$ is well-defined; moreover, it clearly
commutes with $A_r$.}

{From \eqref{trf}, the frame $(V_1,V_2,V_3)$ is obtained (after
rescaling) from the frame $(E,A,H)$ by a complete rotation around
the $E$ axis in time $\pi$. Such a rotation is a generator of
$\pi_1(\mathrm{SO}(3))=\zz/2\zz$. Hence $V$ exists globally around
the $\theta$ circle if and only if the lift of $(E,A,H)$ does not,
i.e., if the spin structure is non-trivial along the fiber
$S^1/\{\pm 1\}$. Otherwise, if the spin structure is trivial along
the fiber $S^1/\{\pm 1\}$, then after time $\pi$ the lift $V$
changes sign.}

{A spinor $[V,\sigma]$ is in the {$m$}-eigenspace of
$i\partial_\theta$ if and only if
\begin{equation}\label{egs}
\sigma(t+\theta)=e^{-im\theta}\sigma(t).
\end{equation}
The resulting spinor should be $\pi$-periodic (since we work on
$\mathrm{PSL}(2,\rz)$, we assumed that $-1\in\Gamma$). We
distinguish two cases:
\begin{itemize}
\item The spin structure is nontrivial along the $S^1/\{\pm 1\}$
fiber. Then $V(\pi)=V(0)$ so we want $\sigma(\pi)=\sigma(0)$. The
eigenspinor equation \eqref{egs} gives $m\in 2\zz$; \item The spin
structure is trivial along the $S^1/\{\pm 1\}$ fiber. Then
$V(\pi)=-V(0)$ so we want $\sigma(\pi)=-\sigma(0)$. The
eigenspinor equation \eqref{egs} gives $m\in 1+2\zz$.
\end{itemize}
In both cases, the $m$-eigenspaces are $2$-dimensional
representation spaces for $c(V_j)$, $j=1,2,3$.}

Denote by $A_{r,m}$ the action of $A_r$ on the $m$-eigenspace of
$i\partial_\theta$. We get a $b$-operator $A_{r,m}$ (in the sense
of Melrose) in $L^2([0,\epsilon),\cz^2,\nu^{-1}d\nu)$
\[A_{r,m}=2c(V_1) {\nu}\partial_{\nu}
-c\left(\frac{V_3}{r}-V_2\right)im {-}\frac r2.\] The $b$-normal
operator of $A_{r,m}$ is obtained by replacing
${\nu}\partial_{\nu}$ with $is$ where $s$ is a complex parameter.
One knows from the general theory of $b$-operators \cite{Mel93}
that the following statements are equivalent:
\begin{itemize}
\item $\lambda$ does not belong to the essential spectrum of
$A_{r,m}$, i.e., $A_{r,m}-\lambda$ is Fredholm; \item
$\Nn(A_{r,m})(s)-\lambda$ is invertible for all $s\in\rz$.
\end{itemize}
We use now the representation
\begin{align*}c(V_1)=\begin{bmatrix}i&0\\0&-i\end{bmatrix},&&
c\left(\frac{V_3}{r}-V_2\right)=\begin{bmatrix}
0&1+r^{-2}\\-1&0\end{bmatrix}
\end{align*}
so that
\[\Nn(A_{r,k})(s)=\begin{bmatrix}-2s&-im(1+r^{-2})\\im&2s\end{bmatrix}
{-}\frac r2.\] An easy computation shows that
$\Nn(A_{r,k})(s)-\lambda$ is invertible for all $s\in\rz$ exactly
for
\[\lambda\in \left(-\frac r2-|m|(1+r^2)^{\frac12}r^{-1},\
-\frac r2+|m|(1+r^2)^{\frac12}r^{-1}\right).\] Thus the
essential spectrum of $D_r$ is the superposition of the
complements of these intervals for each $k$ and for each cusp
$c_p$ with trivial spin structure.

In the sequel, we will assume that the spin structure is trivial
along the $S^1/\{\pm 1\}$ fiber, so the essential spectrum does
not touch $0$ for small $r>0$. {In this case an alternate proof of
Theorem \ref{t:cont-spec} {will follow from a computation using
harmonic analysis over $G$ (see Section \ref{sec-FT})}.}

%%%%%%%%%%%%%%%%%%%%%%%%%%%%%%%%%%%%%%%%%%%%%%%%%%%%%%%%%%%%%%%%%%%%
\section{Regularized trace and Geometric side}
\label{sec-regtracandgeomside}
%%%%%%%%%%%%%%%%%%%%%%%%%%%%%%%%%%%%%%%%%%%%%%%%%%%%%%%%%%%%%%%%%%%%

In this section, we study the relation of certain regularized
trace of the odd heat operator of $\Dd_r$ with the geometric side
of the Selberg trace formula.

To use the harmonic analysis over $G$, we need to fix the Haar
measures over $G$ and its subgroups. First the parametrizations in
\eqref{e:param} for $A_0,N_0$ carry the Lebesgue measure $du, dn$
from $\R$ to $A_0, N_0$.  Now we fix Haar measures on $K$ by
$\mathrm{vol}(K/Z)=1$ and on $G$ by
\begin{equation}\label{e:Haar}
\int_G f(g)\, dg= \int_{N_0}\int_{A_0}\int_K f(na_uk) e^{-u} \,
dk\, {du} \, dn
\end{equation}
for $f_0\in C_0(G)$ and $a_u=\mathrm{diag}(e^{\frac u2}, e^{-\frac
u2})$. For $a_{P,u}:=k_P^{-1} a_u k_P\in A=k_P^{-1} A_0 k_P$, we
put
\[
H_P(g)=u \qquad \text{for}\quad g\in N a_{P,u} K.
\]

The Iwasawa decomposition $\Hh\cong G/K \cong NA$ provides a
parametrization of the geodesics $nA\cdot i\subset \Hh$ to the
infinity. The parameter value is given by the function $H_P$ whose
potential {curves} are $N$-orbits (horocycles) on $\Hh$. However,
this parametrization is not adapted to $\Gamma$. To rectify this,
we replace $k_P$ by $g_P=a_{u_P} k_P$ where
$e^{-u_P}=\mathrm{vol}(\Gamma_N \backslash N)$ where
$\Gamma_N:=\Gamma\cap N$. For the new parameter
\[
H_P(g)+u_P=H_{P_0}(g_P g),
\]
then the value $0$ of this new parameter corresponds to the
horocycle whose projection on $\Gamma\backslash \Hh$ has
length $1$.

For $\phi\in \mathsf{H}:=L^2(Z\backslash K)=\langle e^{im\theta}\,
|\, m\in 2\mathbb{Z}\rangle$ and $s\in\mathbb{C}$, we extend
$\phi$ to $G$ by
\[
\phi_s(na_uk)=e^{su}\phi(k) \qquad\text{for} \quad n\in N_0, k\in
K.
\]
These functions constitute the Hilbert space $\mathsf{H}_s\cong
\mathsf{H}$ in which the representation $\pi_s$ induced from the
parabolic subgroup $P_0=N_0A_0Z$ acts as
\[
(\pi_s(g)\phi_s)(x)=\phi_s(xg).
\]

From now on, we assume that $\mathfrak{P}=\{P_1, \ldots
,P_{\kappa}\}$ is {a} set of representatives for
$\Gamma$-conjugacy classes of the cuspidal parabolic subgroups and
that the spin structure over $c_{P_i}$ for $1\leq i\leq \kappa^\t$
is trivial. We also assume that the representation $\rho$ {maps
the generator $k\in\pi_1(X)$ to the identity $1$, thus we consider
only spin structures which are trivial along the $S^1/\{\pm
1\}$-fiber.} For the representation space $V\cong\mathbb{C}^2$ of
$\chi=\rho\oplus \rho$, we let $V^P$ be the invariant subspace of
$V$ under the action $\chi|_{\Gamma_P}$. Then
\[
V^{P_i}=\begin{cases} V \quad \quad \text{if} \quad 1\leq i \leq \kappa^\t\\
                      \{0\} \ \quad \text{if} \quad \kappa^\t+1\leq
                      i\leq \kappa \end{cases}.
                      \]
For a cuspidal parabolic subgroup $P$, $s\in\mathbb{C}$ with
$\Re(s)>1$ and $\phi\in\mathsf{H}\otimes V^P$, the Eisenstein
series $E(P,\phi,s)$ is defined by
\[
E(P,\phi,s)(g):=\sum_{\gamma\in \Gamma/\Gamma_P} \chi
(\gamma){\phi_s}(g_P\gamma^{-1} g).
\]
Note that there is no Eisenstein series attached to $P_i$ if
$\kappa^\t+1\leq i\leq \kappa$. The Eisenstein series
$E(P,\phi,s)$ converges absolutely and locally uniformly for
$\Re(s)
>1$ and has the meromorphic extension over $\mathbb{C}$. In particular,
$E(P,\phi,s)$ is an \emph{automorphic form}, that is,
\[
E(P,\phi,s)(\gamma g)=\chi(\gamma) E(P,\phi,s)(g) \qquad
\text{for} \quad \gamma\in\Gamma, g\in G.
\]
For $\phi\in H\otimes V^{\mathrm{cst}}$ with
$V^{\mathrm{cst}}:=\oplus_{P\in\mathfrak{P}} V^P$, we define
\begin{equation}\label{e:Eisenstein}
E(\phi,s)=\sum_{P\in\mathfrak{P}}E(P,\mathrm{pr}^P\phi,s)\end{equation}
where $\mathrm{pr}^P$ denotes the orthogonal projection onto
$V^P$, and
\[
E^{\mathrm{cst}}(\phi,s)(g)=
\big(E^P(\phi,s)(g_P^{-1}g)\big)_{P\in\mathfrak{P}}.\]
Here, the \emph{constant term} of $E^P(\phi,s)$ is defined by
\[
E^P(\phi,s)(g):=\mathrm{vol}(\Gamma_N\backslash
N)^{-1}\int_{\Gamma_N\backslash N} \mathrm{pr}^P\, E(\phi,s)(ng)\,
dn
\]
{for $N=N_P$.} Then we have
\[
E^\cst(\phi,s)=\phi_s+ \big(C(s)\phi)_{1-s}
\]
where $C(s)$ is the \emph{scattering operator} acting on
$\mathsf{H}\otimes V^{\mathrm{cst}}$.

Now let us describe the spectral decomposition of $\sLchi$,
\[
\sLchi=\cusLchi\oplus \resLchi \oplus \ctLchi.
\]
Here $\cusLchi$ is the space of the cusp forms in
$L^2(\Gamma\backslash G, \chi)$, and decomposes into a Hilbert
direct sum of closed irreducible $G$-invariant subspaces with
finite multiplicities. The residual part $\resLchi$ is the direct
sum of the constants and of finitely many copies of the
complementary series representation of $G$ such that some
Eisenstein series has a pole at $s\in (\frac12,1)$. These two
spaces constitute the discrete part $\disLchi$. The continuous
part $\ctLchi$ is isometric to
\[
\Big\{\, \Phi\in
L^2(\frac12+i\mathbb{R},\frac{d\tau}{4\pi})\hat{\otimes}
\mathsf{H}\otimes V^\cst\, | \
\Phi(\frac12-i\tau)=C(\frac12+i\tau)\Phi(\frac12+i\tau)\, \Big\}
\]
by
\begin{equation}\label{e:intert}
\mathcal{I}^{\mathrm{ct}}\Phi=\frac{1}{4\pi}\int^\infty_{-\infty}
E(\Phi, \frac12+i\tau) \, d\tau
\end{equation}
where $E(\Phi,\frac12+i\tau)$ is defined as in
\eqref{e:Eisenstein} with $\phi=\Phi$ and $s=\frac12+i\tau$. For
$f\in L^1(G)$, we define a representation on $\sLchi$ by
\[
\pi(f):=\int_G f(g)\pi(g)\, dg
\]
where $\pi$ is the right translation action given by
$\big(\pi(g)\phi\big)(x)=\phi(xg)$ for $\phi\in\sLchi$. We put
$\pi^{\dis}(f)=\pr^\dis\circ \pi(f)$,
$\pi^{\ct}(f)=\pr^{\ct}\circ\pi(f)$ where $\pr^{\dis},\pr^{\ct}$
denote the orthogonal projections onto $\disLchi,\ctLchi$
respectively. In particular, $\pi^\ct(f)$ intertwines with
$\pi^\cst_{\frac12+i\tau}(f):=\pi_{\frac12+i\tau}(f)\otimes
\mathrm{Id}_{V^\cst}$ by $\mathcal{I}^\ct$ given in
{{\eqref{e:intert}}.

Since $\widetilde{D}_r$ is a left invariant differential operator,
there is a function ${\hat{f}_{t,r}}\in C^\infty(G,{M(2,\cz)})$
such that
\[
{\hat{f}}_{t,r}(x^{-1}y)= \big(\widetilde{D}_r
e^{-t\widetilde{D}^2_r}\big) (x, y)
 \qquad \text{for}\quad x,y\in G.
\]
By the heat kernel estimates in \cite{Donnelly}, which also holds
for the generalized Laplacian $\widetilde{D}^2_r$ with the
form in \eqref{e:Lap}, we have
\[
|| d_t^i\,d_x^j\, d_y^k\, \hat{f}_{t,r}(x^{-1}y)|| \leq C
t^{-\frac 52-i-j-k}\, \exp\Big(-\frac{d_G^2(x,y)}{4t}\Big)
\]
where $C$ is a positive constant and $d_G$ is the metric over $G$.
(Note that we apply the method in \cite{Donnelly} to a certain
co-compact discrete subgroup $\Gamma'$ in $G$ to obtain the above
estimate.) This estimate implies that
${f_{t,r}}:=\tr(\hat{f}_{t,r})$ lies in the Harish-Chandra
$L^1$-Schwartz space $\mathcal{C}^1(G) (\subset L^1(G))$ defined
by
\begin{multline*}
\Cc^1(G)= \Big\{ f\in C^\infty(G)\ \big|\
|f(D_1k_{\theta_1}a_uk_{\theta_2}D_2)| \\ \leq
Ce^{-|u|}(1+|u|+|\theta_1+\theta_2|)^{-n}, \quad \forall
n\in\mathbb{N}, D_1,D_2\in \mathfrak{g} \Big\}
\end{multline*}
where $f(D_1k_{\theta_1}a_uk_{\theta_2}D_2)$ denotes the
convolution
$D_1*\delta_{k_{\theta_1}}*\delta_{a_u}*\delta_{k_{\theta_2}}*D_2$
evaluated on $f$. Let us put
\begin{align*}
K(t,x,y):= \sum_{\gamma\in Z\backslash\Gamma}
\hat{f}_{t,r}(x^{-1}\gamma y)\,\chi(\gamma)=\sum_{\gamma\in
Z\backslash\Gamma} \widetilde{D}_r e^{-t\widetilde{D}^2_r} (x,
\gamma y)\, \chi(\gamma) && \text{for $x,y\in G$},
\end{align*}
which is absolutely uniformly convergent on compact sets in $G$.
For a $\Gamma$-cuspidal parabolic subgroup $P_j=P=NAZ$, we define
the constant term of $K(t,x,y)$ along $P$ as follows,
\[
K_P(t,x,y)= \mathrm{vol}(\Gamma_N\backslash N)^{-1}
\int_{\Gamma_N\backslash N} \sum _{\gamma\in Z\backslash \Gamma_P}
\hat{f}_{t,r}(x^{-1} \gamma n y)\, \chi(\gamma)\, \pr^P \, dn .
\]
For $u\in\mathbb{R}$, let $\alpha_{P}(u)$ be the characteristic
function of $\{ x\in G \, | \, H_P(x)+u_P> u\}$, which projects on
certain region $\mathcal{C}_{P,u}\subset \Gamma\backslash G$ for a
large $u$. Then the truncation of $K(x,x)$ is defined by
\[
\Lambda_u K(t,x,x):= K(t,x,x) - \sum_{P\in\mathfrak{P}}
\alpha_{P}(u) K_{P}(t,x,x),
\]
which is an automorphic form over $\Gamma\backslash G$.

\begin{proposition}\emph{(Maass-Selberg Relation)} \label{p:Maass-Sel}
For $u\gg 0$, we have
\begin{multline}\label{e:Maass-Sel}
\int_{\Gamma\backslash G} \tr \big( \Lambda_u K(t,x,x) \big) \, dx
\,\\=\,   u\frac{1}{2\pi}
\int^\infty_{-\infty}\Tr\big(\pi^\cst_{\frac12+i\tau}(f_{t,r})\big)
\, d\tau + \Tr\big( \pi^{\dis}(f_{t,r})\big) + \frac{1}{4}
\Tr\big(
C({\frac{{1}}{{2}}})\pi^{\mathrm{cst}}_{\frac12}(f_{t,r})\big)  \\
 -\frac{1}{4\pi}\int^\infty_{-\infty}\,
\Tr \big( C'(\frac12-i\tau)C(\frac12+i\tau)
\pi_{\frac12+i\tau}^{\mathrm{cst}} (f_{t,r})\big) \, d\tau.
\end{multline}
\end{proposition}
\begin{proof} For {a} test function with compact support and
$K$-finite condition, we can prove this proposition just following
the argument in p.58--60 of \cite{Hoff94}. Then this can be
generalized easily to our test function $f_{t,r}$ as in proof of
the theorem 25 of \cite{Hoff94}. The finiteness of the integrand
of the integrals on the right hand side follows from Theorem
\ref{t:cont-spec}.
\end{proof}

From Proposition \ref{p:Maass-Sel}, one can see that the first
term on the right side of \eqref{e:Maass-Sel} is blowing up as
$u\to\infty$. Hence it is natural to remove this term in the
following definition,
\begin{align}\label{e:reg-tr}
\Tr\big(\Dd_r e^{-t\Dd_r^2}\big):=& \Tr\big(
\pi^{\dis}(f_{t,r})\big) + \frac{1}{4} \Tr\big(
C(\frac12)\pi^{\mathrm{cst}}_{\frac12}(f_{t,r})\big) \\
&\qquad\qquad\ -\frac{1}{4\pi}\int^\infty_{-\infty}\, \Tr \big(
C'(\frac12-i\tau)C(\frac12+i\tau)
\pi_{\frac12+i\tau}^{\mathrm{cst}} (f_{t,r})\big)\, d\tau.\notag
\end{align}
This regularized trace is the essentially same as the
$b$-trace of Melrose \cite{Mel93}, and is related with the
geometric side of the Selberg trace formula as we will see in
Proposition \ref{p:Selberg}. Denote
\[
h(\tau)=\Theta_{\frac12+i\tau}(f), \quad h(n)=\Theta_{n}(f)
\] where
$\Theta_{\frac12+i\tau}(f):=\Tr\big(\pi_{\frac12+i\tau}(f)\big)$
for a principal series representation $\pi_{\frac12+i\tau}$, and
$\Theta_n(f):=\Tr\big(\pi_n(f)\big)$ for a discrete series
representation $\pi_n$. An operator $J(s)$ over
$\mathsf{H}=L^2(Z\backslash K)$ is defined by
\begin{equation}\label{e:J}
J(s)\phi_m= \frac{1}{\sqrt{\pi}}\cdot
\frac{\Gamma(s)\Gamma(s-\frac12)}{\Gamma(s+\frac m2)\Gamma(s-\frac
m2)}\, \phi_m
\end{equation}
for the basis $\phi_m(k_\theta)=e^{im\theta}\in
\mathsf{H}$ where $k_\theta= \begin{pmatrix} \cos\theta & \sin\theta\\
-\sin\theta &\cos\theta \end{pmatrix}$. The following proposition
follows from Theorem 13 and Lemma 24 in \cite{Hoff94}.

\begin{proposition} \label{p:Selberg}\emph{(Selberg Trace
Formula)}
\begin{align}\label{e:Selberg}
\Tr\big(\Dd_r e^{-t\Dd_r^2}\big) =&
\frac{\mathrm{vol}(\Gamma\backslash G)}{2\pi} \Big(
\int^\infty_{-\infty} {\tau\tanh (\pi \tau)}\, h_{t,r}(\tau)
\, d\tau + \sum_{n\equiv 0 (\mathrm{mod} 2)} (|n|-1) h_{t,r}(n) \Big)\\
&+ \sum_{[\gamma]\in Z\backslash \Gamma_{\mathrm{hyp}}}
\frac{\tr\big(\chi(\gamma)\big)\, u_\gamma}{4\pi[\Gamma_\gamma:Z]
\sinh
\frac{u_\gamma}{2}} \int^\infty_{-\infty} \cos(u_\gamma \tau)\, h_{t,r}(\tau)\, d\tau\notag\\
&-2\kappa^{\t} \Big(\frac{1}{2\pi}\int^\infty_{-\infty}
\psi(1+2i\tau)h_{t,r}(\tau)\, d\tau +\frac12\sum_{n\equiv 0
(\mathrm{mod} 2)} h_{t,r}(n)
\Big)\notag\\
&+2(\kappa-\kappa^\mathrm{t}) \frac{\log
2}{2\pi}\int^\infty_{-\infty} h_{t,r}(\tau) \, d\tau\\
 &+\frac{\kappa^\t}{2}
h_{t,r}(0)-\frac1{4\pi}\mathrm{p.v.}\int^\infty_{-\infty}
\Tr\big(J(\frac12+i\tau)^{-1} J'(\frac12+i\tau)
\pi_{\frac12+i\tau}(f_{t,r}) \big) \, d\tau\notag
\end{align}
where $h_{t,r}(\tau),h_{t,r}(n)$ are defined for $f_{t,r}$,
the sum $\sum_{[\gamma]\in \Gamma_{\mathrm{hyp}}}$ is given
over the $\Gamma$-conjugacy class of hyperbolic elements $\gamma$
conjugate to $a_{u_\gamma}$, and
$\psi(z)=\Gamma'(z)\Gamma(z)^{-1}$.
\end{proposition}

%%%%%%%%%%%%%%%%%%%%%%%%%%%%%%%%%%%%%%%%%%%%%%%%%%%%%%%%%%%%%%%%%
\section{Fourier transforms $h_{t,r}(\tau),h_{t,r}(n)$}
\label{sec-FT}
%%%%%%%%%%%%%%%%%%%%%%%%%%%%%%%%%%%%%%%%%%%%%%%%%%%%%%%%%%%%%%%%%

In this section, we compute $h_{t,r}(\tau),h_{t,r}(n)$ which are
needed to analyze the right hand side of the Selberg trace
formula.

First, let us consider $h_{t,r}(\tau)$. {For this, recall
\begin{equation}\label{e:htr-old}
h_{t,r}(\tau)=\Tr\big(\pi_{\frac12+i\tau}(f_{t,r})\big)=\sum_{n=1}^\infty\int_G
  f_{t,r}(g)  \left(\pi_{\frac12+i\tau}(g)
\xi_n,\xi_n \right) \, dg
\end{equation}
where $\{\xi_n\}_{n=1}^\infty$ is the orthonormal basis of the
representation space of $\pi_{\frac12+i\tau}$, which is given by
the union of the following spaces indexed by $m\in\mathbb{Z}$ for
$s=\frac12+i\tau$,
\[
\mathsf{H}(s,m):=\big\{\, \phi_{s}\in \mathsf{H}_{s}\, \big|\,
\phi_{s}(na_uk_\theta)=e^{su}e^{im\theta}\ \ \text{for}\
na_uk_\theta\in N_0A_0K \, \big\}.
\]
Since $\widetilde{\Dd}_r$ is $Z$-invariant, the Fourier transform
$h_{t,r}(\tau)$ is nontrivial only if $m$ is an even number.}
Recalling
\[
{\widetilde{\Dd}_r} = \Big(\frac{2-r^2}{2r}\Big)+
\begin{pmatrix} r^{-1}Z & 2X_-
\\ -2X_+ & -r^{-1}Z \end{pmatrix}
\]
the problem is again reduced to the following lemma, which can be
obtained applying the equalities in \eqref{e:diff}.

\begin{lemma}\label{l:action} We have
\begin{equation*}
Z f= m f, \quad X_{\pm} f = -\frac{i}{2}\big(m\pm 2s \big)
e^{\pm2i\theta} f \qquad \text{for $f\in \mathsf{H}(s,m)$}.
\end{equation*}
\end{lemma}

The second equality in Lemma \ref{l:action} implies that $X_\pm$
maps $\mathsf{H}(s,m)$ to $\mathsf{H}(s,m\pm 2)$.  From these
facts,
\begin{equation*}
{\widetilde{\Dd}_r} \begin{pmatrix} \phi_{\tau,m-2} \\
\phi_{\tau,m}\end{pmatrix}
=\begin{pmatrix} r^{-1}(m-2)+ 2^{-1}\ell & -{i}(m -1 -2i\tau) \\
i(m-1 +2i\tau) & -r^{-1}m + 2^{-1}\ell \end{pmatrix}
\begin{pmatrix} \phi_{\tau,m-2} \\ \phi_{\tau,m}\end{pmatrix}
\end{equation*}
where $\ell=\frac{2-r^2}{r}$ and $\phi_{\tau,m-2}\in
\mathsf{H}(\frac12+i\tau,m-2),\phi_{\tau,m}\in\mathsf{H}(\frac12+i\tau,m)$.
Hence the action of ${\widetilde{\Dd}_r}$ on
$\mathsf{H}(\frac12+i\tau,m-2)\oplus \mathsf{H}(\frac12+i\tau,m)$
is given by the roots of
\begin{equation}\label{e:eigen1}
\lambda^2+r\lambda + \frac{r^2}{4}-\frac{(m-1)^2}{r^2}=
(m-1)^2+4\tau^2,
\end{equation}
that is,
\[
\lambda_\pm(\tau,m)= -\frac{r}{2} \pm
\big((m-1)^2(1+r^{-2})+4\tau^2\big)^{1/2} \qquad \text{for $m\in
2\,\Z, \tau\in \R^+$}.
\]
Therefore we have

\begin{lemma}\label{l:cont-htr}
\[
h_{t,r}(\tau)=\Theta_{\frac12+i\tau}(f_{t,r})=\sum_{m\in
2\mathbb{Z}} \big(\lambda_+(\tau,m) e^{-t\lambda_+(\tau,m)^2}
+\lambda_-(\tau,m)e^{-t\lambda_-(\tau,m)^2}\big).
\]
\end{lemma}

Just repeating the above computation applied to the
Eisenstein series $E(P,\phi_m,s)$, we can also prove Theorem
\ref{t:cont-spec}.

Next we compute $h_{t,r}(n)$ for the discrete series
representation $\pi_n$. We review the discrete series
representations of $G=\mathrm{SL}(2,\R)$. For this it is more
convenient to use the Lie group $\mathrm{SU}(1,1)$ which is
conjugate to $\mathrm{SL}(2,\R)$ within
$\mathrm{SL}(2,\mathbb{C})$:
\begin{align*}
\begin{pmatrix} 1 & i \\ i & 1 \end{pmatrix} \mathrm{SU}(1,1) \begin{pmatrix} 1 & i \\ i & 1
\end{pmatrix}^{-1}
= \mathrm{SL}(2,\R).
\end{align*}
Here
\begin{equation*}
\mathrm{SU}(1,1) = \left\{ \begin{pmatrix} \alpha & \beta \\
\overline{\beta} & \overline{\alpha} \end{pmatrix} \Big\vert \ \
\vert \alpha\vert^2 - \vert \beta \vert^2 =1 \right\}.
\end{equation*}
Then the holomorphic discrete series $\pi_n$ ($n\in\mathbb{N}$) as
a representation of $\mathrm{SU}(1,1)$ acts on analytic functions
on the disc by
\begin{equation*}
\pi_n \begin{pmatrix} \alpha & \beta \\
\overline{\beta} & \overline{\alpha} \end{pmatrix} f(z) = (-\beta
z +\overline{\alpha})^{-n} f\left( \frac{\alpha z
-\overline{\beta}}{-\beta z +\overline{\alpha}}\right),
\end{equation*}
and the norm, except for a constant factor, is given by
\begin{equation*}
\vert\vert f \vert\vert = \begin{cases} &\int_{\vert z\vert < 1}
\left\lvert f(z) \right\rvert^2 (1 -\vert z \vert^2)^{n-2} \, dz
\qquad \text{for $n\ge 2$}\\
& \text{sup}_{0 \le r < 1} \int^{2\pi}_0 \left\lvert f(r
e^{i\theta}) \right\rvert ^2 \, d\theta \qquad \quad \text{for
$n=1$}
\end{cases}.
\end{equation*}
The anti-holomorphic discrete series $\pi_n$ ($n\in-\mathbb{N}$)
as a representation of $\mathrm{SU}(1,1)$ acts on analytic
functions on the disc by
\begin{equation*}
\pi_n \begin{pmatrix} \alpha & \beta \\
\overline{\beta} & \overline{\alpha} \end{pmatrix} f(z) =
(-\overline{\beta} z +{\alpha})^{-n} f\left(
\frac{\overline{\alpha} z -{\beta}}{-\overline{\beta} z
+{\alpha}}\right)
\end{equation*}
with the same norm.

\begin{lemma}\label{p:dis-act}For the basis $\{z^N\}_{N\in\{0\}\cup \mathbb{N}}$
of the space of analytic functions on the disc, we have
\begin{multline*}
X_+ z^N =\, (N+n) z^{N+1}, \quad X_- z^N =-N z^{N-1}, \\  Z\, z^N
=
 (2N+n)\, z^N \qquad \text{by the action of $\pi_n$, $n\in\mathbb{N}$},
\end{multline*}
\begin{multline*}
 X_+ z^N = -N z^{N+1}, \quad
X_-z^N = (N+n) z^{N-1}, \\ Z\, z^N = - (2N+n) z^N \qquad \text{by
the action of $\pi_n$, $n\in-\mathbb{N}$}.
\end{multline*}
\end{lemma}

\begin{proof}
By elementary computations, we can see that the subgroups
generating $K,H,A$ are transformed as follows:
\begin{align*}
\begin{pmatrix} \cos\theta & \sin\theta \\ -\sin\theta &
\cos\theta \end{pmatrix} \qquad &\longrightarrow \qquad
\begin{pmatrix} e^{i\theta} & 0 \\ 0 & e^{-i\theta} \end{pmatrix}
\\
\begin{pmatrix} e^{t} & 0 \\ 0 & e^{-t} \end{pmatrix}\qquad \qquad
&\longrightarrow \qquad \begin{pmatrix} \cosh t & i\sinh t \\
-i\sinh t & \cosh t \end{pmatrix} \\
\begin{pmatrix} \cosh t & \sinh t \\ \sinh t & \cosh
t\end{pmatrix} \qquad &\longrightarrow \qquad
\begin{pmatrix} \cosh t & \sinh t \\ \sinh t & \cosh
t\end{pmatrix}
\end{align*}
where the matrices on the right side denote elements in
$\mathrm{SU}(1,1)$. To see the action of $K$ under $\pi_n$ for
$n\in\mathbb{N}$, let us consider
\begin{equation}\label{e:weight}
\pi_n \begin{pmatrix} e^{i\theta} & 0 \\
0 & e^{-i\theta} \end{pmatrix} z^N = e^{(2N+n)i\theta} z^N,
\end{equation}
which implies
\begin{equation*}
Z\, z^N = -i\, \frac{d}{d\theta}\Big\vert_{\theta=0}
\pi_n \begin{pmatrix} e^{i\theta} & 0 \\
0 & e^{-i\theta} \end{pmatrix} z^N = (2N+n)\, z^N .
\end{equation*}
In a similar way, we can show that the action $\pi_n$ for
$n\in\mathbb{N}$ by $H,A$ are given by
\begin{align*}
H\, z^N =i(N+n) z^{N+1} + i N z^{N-1},\qquad A\, z^N = (N+n)
z^{N+1} - N z^{N-1}.
\end{align*}
These imply the equalities for $\pi_n$ for $n\in\mathbb{N}$. The
case for $\pi_n$ for $n\in-\mathbb{N}$ can be obtained by taking
the complex conjugates of equalities for the action of $\pi_n$ for
$n\in\mathbb{N}$.
\end{proof}

Now we consider the action of ${\widetilde{\Dd}_r}$ under $\pi_n$
for $n\in\mathbb{N}$. From \eqref{e:weight}, we can see that $z^N$
are the $K$-type vectors of weight $m$ if $m=2N+n$. By Proposition
\ref{p:dis-act}, over $(\alpha,\beta)$ for $K$-type $(m-2),m$
vectors $\alpha,\beta$ in the representation space of $\pi_n$ the
Dirac operator ${\widetilde{\Dd}_r}$ acts  by
\begin{equation*}
\biggl(\frac{2-r^2}{2r}\biggr)\mathrm{Id}+ \begin{pmatrix}
r^{-1}(m-2) &  n-m
\\ -(n+m-2) & -r^{-1}m \end{pmatrix}
\end{equation*}
noting $N= (m-n)/2$. We have two cases: First, if $\beta$ is not
the minimal $K$-type for $\pi_n$, that is, $m\gneq n$, as in the
derivation of \eqref{e:eigen1}, we can obtain the corresponding
eigenvalue equation
\[
\lambda^2+r\lambda + \frac{r^2}{4}-\frac{(m-1)^2}{r^2}=
(m-1)^2-(n-1)^2.
\]
Hence
\[
\lambda_\pm(n,m)= -\frac{r}{2} \pm
\big((m-1)^2(1+r^{-2})-(n-1)^2\big)^{1/2} \qquad \text{for
$m=n+2,n+4,\ldots$}.
\]
Second, if $\beta$ is the minimal $K$-type for $\pi_n$, that is,
$K$-type $m=n$ vector, then $\alpha$ is just trivial. Hence, the
eigenvalue is given by
\[
\la(n)=-\frac{r}{2}+\frac{1-n}{r}.
\]
Repeating the same procedure as in the case of $\pi_n$ for
$n\in\mathbb{N}$, we can obtain the same eigenvalues
$\lambda_\pm(n,m)$, $\la(n)$ for $\pi_n$ for $n\in-\mathbb{N}$. In
our case, $n$ should be an even number {since} $\Gamma$ contains
$-\mathrm{Id}$. Combining all the facts derived in the above, we
have

\begin{lemma}\label{l:dis-htr} For $n\in2\mathbb{N}, h(n)=\Theta_{n}(f)$, we have
\begin{multline*}
h_{t,r}(n)=h_{t,r}(-n)=\Big(\lambda(n)e^{-t\lambda(n)^2}\\+\sum_{m\in
n+2\mathbb{N}} \big(\lambda_+(n,m) e^{-t\lambda_+(n,m)^2}
+\lambda_-(n,m)e^{-t\lambda_-(n,m)^2}\big) \Big)
\end{multline*}
where $\la(n)=-\frac{r}{2}+\frac{1-n}{r}, \ \lambda_\pm(n,m)=
-\frac{r}{2} \pm \big((m-1)^2(1+r^{-2})-(n-1)^2\big)^{1/2}$.
\end{lemma}

%%%%%%%%%%%%%%%%%%%%%%%%%%%%%%%%%%%%%%%%%%%%%%%%%%%%%%%%%%%%%%%%%%%
\section{Eta function of $\Dd_r$: Principal series part}
\label{sec-etafnprincipalseriespart}
%%%%%%%%%%%%%%%%%%%%%%%%%%%%%%%%%%%%%%%%%%%%%%%%%%%%%%%%%%%%%%%%%%%

Now we study the eta function defined by
\[
\eta({\Dd_r}, s):=\frac{1}{\Gamma(\frac{s+1}{2})}\int^\infty_0
t^{\frac{s-1}{2}} \Tr\big(\Dd_r e^{-t\Dd_r^2}\big)\, dt
\]
for $\Re(s)\gg 0$ and $r$ near $0$. First let us recall that the
bottom of each branch of continuous spectrum of $\Dd_r$ goes to
$\infty$ as $r\to 0$ by Theorem \ref{t:cont-spec}. Hence for a
small $r>0$, $\Tr\big(\Dd_r e^{-t\Dd_r^2}\big)$ decays
exponentially as $t\to\infty$. To analyze $\Tr\big(\Dd_r
e^{-t\Dd_r^2}\big)$ in more detail, we apply Proposition
\ref{p:Selberg} which relates $\Tr\big(\Dd_r e^{-t\Dd_r^2}\big)$
with the geometric side. Note that this geometric side can be
decomposed into two parts:
\begin{align*}
\Tr_\p\big(\Dd_r e^{-t\Dd_r^2}\big)=&\Tr\big(\Dd_r
e^{-t\Dd_r^2}\big)-\Tr_\d\big(\Dd_r e^{-t\Dd_r^2}\big),\\
\Tr_\d\big(\Dd_r e^{-t\Dd_r^2}\big) =&
\frac{\mathrm{vol}(\Gamma\backslash G)}{2\pi}  \sum_{n\equiv 0
(\mathrm{mod} 2)} (|n|-1) h_{t,r}(n)  -\kappa^{\t} \sum_{n\equiv 0
(\mathrm{mod} 2)} h_{t,r}(n)
\end{align*}
and accordingly we also decompose the eta function $\eta(\Dd_r,s)$
into
\[
\eta(\Dd_r,s)=\eta_\p(\Dd_r,s)+\eta_\d(\Dd_r,s).
\]
The principal part of the eta function $\eta_\p(\Dd_r,s)$ is
studied in this section and the other part $\eta_\d(\Dd_r,s)$ will
be considered {in} the next section.

We start with the following lemma.

\begin{lemma}\label{l:h-tau} Putting
$I(m,r,\tau)=\big((m-1)^2(1+r^{-2})+4\tau^2\big)^{\frac12}$,
\begin{multline*}
h_{t,r}(\tau)= \exp\big(-\frac{r^2}{4}t-4\tau^2 t\big)\\
\cdot \sum_{m\in2\mathbb{Z}}e^{-(m-1)^2(1+r^{-2})t}\,
\Big(\sum_{k=0}^\infty\big(-r\frac{(rt)^{2k}}{(2k)!}+
2\frac{(rt)^{2k-1}}{(2k-1)!}\big)\, I(m,r,\tau)^{2k} \Big)
\end{multline*}
where the term $\frac{(rt)^{2k-1}}{(2k-1)!}$ for $k=0$ vanishes
and for $t\in[0,1]$ and $r\in(0,1]$ the following estimate holds,
\begin{multline}\label{e:htr-est}
|h_{t,r}(\tau)|\, \leq\, 2r\exp\big(-\frac{r^2}{4}t-4\tau^2
t\big)\\
\cdot \sum_{m\in2\mathbb{Z}}e^{-(m-1)^2(1+r^{-2})t}\,
\Big(1+I(m,r,\tau)^2+e^{ I(m,r,\tau)^2rt}\Big).
\end{multline}
\end{lemma}

\begin{proof}
We can rewrite $h_{t,r}(\tau)$ as follows,
\begin{multline}\label{e:htr}
h_{t,r}(\tau)=\exp\big(-\frac{r^2}{4}t-4\tau^2 t\big)
\sum_{m\in2\mathbb{Z}}e^{-(m-1)^2(1+r^{-2})t}\\
\cdot \Big(\, -\frac{r}{2} \big(
e^{I(m,r,\tau)rt}+e^{-I(m,r,\tau)rt}\big)
+I(m,r,\tau)\big(e^{I(m,r,\tau)rt}-e^{-I(m,r,\tau)rt}\big)\,
\Big).
\end{multline}
Now the Taylor expansion of $e^{I(m,r,\tau)rt}\pm
e^{-I(m,r,\tau)rt}$ gives us the claimed form of the first
equality. To prove the second estimate, we note that
\begin{multline*}
\sum_{k=0}^\infty\big(-r\frac{(rt)^{2k}}{(2k)!}+
2\frac{(rt)^{2k-1}}{(2k-1)!}\big)\, I(m,r,\tau)^{2k}\\
= -r+ rt\big(2-\frac{r^2t}{2}\big) I(m,r,\tau)^2+
\sum_{k=2}^{\infty}
\frac{(rt)^{2k-1}}{(2k-1)!}\big(2-\frac{r^2t}{2k}\big)
I(m,r,\tau)^{2k}.
\end{multline*}
For $t\in[0,1]$ and $r\in(0,1]$, observe that
\begin{align*}
\sum_{k=2}^{\infty}
\frac{(rt)^{2k-1}}{(2k-1)!}\big(2-\frac{r^2t}{2k}\big)
I(m,r,\tau)^{2k} \, \leq  \, 2r \sum_{k=0}^{\infty}
\frac{(rt)^{k}}{k!} I(m,r,\tau)^{2k},
\end{align*}
from which it is easy to derive the estimate.
\end{proof}

Our first task in this section is to get the asymptotic expansion
of $\Tr_\p\big(\Dd_r e^{-t\Dd_r^2}\big)$ as $t\to 0$. By Lemma
\ref{l:h-tau}, we can rewrite the first part of $\Tr_\p\big(\Dd_r
e^{-t\Dd_r^2}\big)$ as follows,
\begin{align}\label{e:iden}
\int^\infty_{-\infty} {\tau\tanh (\pi \tau)}\, h_{t,r}(\tau) \,
d\tau=& \sum_{m\in2\mathbb{Z}}
\exp\big(-\frac{r^2}{4}t-(m-1)^2(1+r^{-2})t\big)\,\\
&\cdot\sum_{k=0}^\infty \sum_{k=p+q}\big(a_{k,p,q}(r)
t^{2k}+b_{k,p,q}(r)t^{2k-1}\big)\notag\\
&\cdot(m-1)^{2p}(1+r^{-2})^{p}\, \int^\infty_{-\infty} {\tau\tanh
(\pi \tau)}\, (2\tau)^{2q} e^{-4\tau^2 t}\, d\tau \notag
\end{align}
where $a_{k,p,q}(r),b_{k,p,q}(r)$ $(b_{0,p,q}(r)=0)$ depend
only on $r$ {and are of order $\mathrm{O}(r)$ for small $r>0$}.
The integral in the last line can be handled as follows,
\begin{align*}
\int^\infty_{-\infty} {\tau\tanh (\pi \tau)}\, (2\tau)^{2q}
e^{-4\tau^2 t}\, d\tau =& (-1)^q\partial_t^q \int^\infty_{-\infty}
 {\tau\tanh (\pi \tau)}\, e^{-4\tau^2 t}\, d\tau\\
 = &(-1)^q\partial_t^q \int^\infty_0 \tanh(\pi\sqrt{x}) e^{-4tx}\, dx\\
 =&(-1)^q\partial_t^q \frac{\pi}{8t}\int^\infty_0
 \Big(\sum_{k=0}^\infty\frac{(-4tx)^k}{k!} \Big) \frac{\cosh^{-2}(\pi
 \sqrt{x})}{\sqrt{x}}\, dx.
\end{align*}
Hence we have
\begin{equation}\label{e:asymp1}
\int^\infty_{-\infty} {\tau\tanh (\pi \tau)}\, (2\tau)^{2q}
e^{-4\tau^2 t}\, d\tau \, \sim\, \sum_{k=0}^\infty a_k\,
t^{-q-1+k} \qquad\text{as}\quad t\to 0
\end{equation}
{where $a_k$ are independent of $r$.}
 Now we note the following equalities
for the first and the third factors on the right hand side of \eqref{e:iden},
\begin{align}\label{e:poisson}
&\sum_{m\in2\mathbb{Z}}
(m-1)^{2p}(1+r^{-2})^{p}\,\exp\big(-(m-1)^2(1+r^{-2})t\big)\\
=&(-1)^p\partial^p_t\sum_{m\in\mathbb{Z}}\exp\big(-4(m-\frac12)^2(1+r^{-2})t\big) \notag\\
=&(-1)^p\partial^p_t\sum_{m\in\mathbb{Z}}(-1)^m
\frac{\sqrt{\pi}}{2\sqrt{(1+r^{-2})t}}
\exp\big(-\frac{\pi^2m^2}{4(1+r^{-2})t}\big)\notag\\
=&(-1)^p\partial^p_t \left( \frac{\sqrt{\pi}
r}{2\sqrt{(1+r^{2})t}} + \sum_{m\in\mathbb{Z}-\{0\}}(-1)^m
\frac{\sqrt{\pi}}{2\sqrt{(1+r^{-2})t}}
\exp\big(-\frac{\pi^2m^2}{4(1+r^{-2})t}\big)\right)\notag
\end{align}
where the second equality is the Poisson summation formula.
Note that the terms for nonzero $m$ in the last line of
\eqref{e:poisson} decays exponentially as $t\to 0$, so that small
time asymptotics is given by the first term in the last line of
\eqref{e:poisson}. Therefore we have
\begin{equation}\label{e:asymp2}
\sum_{m\in2\mathbb{Z}}
(m-1)^{2p}(1+r^{-2})^{p}\,\exp\big(-(m-1)^2(1+r^{-2})t\big)\,
\sim\, {a(r)}\, t^{-\frac12-p}  \qquad\text{as}\quad t\to 0
\end{equation}
{where $a(r)$ depends only on $r$ and $\mathrm{O}(r)$ for small
$r>0$}.

 By {\eqref{e:iden} and the asymptotic expansions in} \eqref{e:asymp1}, \eqref{e:asymp2}, {taking care of $r$-dependence of
 their
 coefficients,}
 we can conclude
\begin{equation}\label{e:asymp3}
\int^\infty_{-\infty} {\tau\tanh (\pi \tau)}\, h_{t,r}(\tau) \,
d\tau \, \sim\, \sum_{k=0}^\infty {a_k(r)}\, t^{-\frac32+k}
\qquad\text{as}\quad t\to 0
\end{equation}
{where $a_k(r)$ depends only on $r$ and is of $\mathrm{O}(r^2)$
for small $r>0$.}

 Now for the second part of $\Tr_\p\big(\Dd_r
e^{-t\Dd_r^2}\big)$, we repeat the above process and noting
\begin{multline*}
\int^\infty_{-\infty} \cos(u_\gamma \tau)\,(2\tau)^{2q}
e^{-4\tau^2 t} \, d\tau=(-1)^q \partial^q_t\int^\infty_{-\infty}
\cos(u_\gamma \tau)\, e^{-4\tau^2 t} \, d\tau\\
= (-1)^q \partial^q_t \,\Big(\,\frac{\sqrt{\pi}}{\sqrt{t}}
\exp(-\frac{u_\gamma^2}{4t})\, \Big),
\end{multline*}
we can see that this term does not contribute to the asymptotics
as $t\to 0$.

To deal with the third part of $\Tr_\p\big(\Dd_r
e^{-t\Dd_r^2}\big)$, we recall
\begin{equation}\label{e:digamma}
\psi(1+z) \, \sim \, \log z + \frac{1}{2z}
-\sum^\infty_{k=1}\frac{B_{2k}}{2k} z^{-2k} \qquad \text{as}\quad
z\to\infty
\end{equation}
where $B_{2k}$ is the Bernoulli's number, which implies
\begin{align}\label{e:asymp4}
\int^\infty_{-\infty} \psi(1+2i\tau) &\, (2\tau)^{2q} e^{-4\tau^2
t}\, d\tau \\ &\sim\, b\, t^{-\frac12-q}\log t\, + c\, +
\sum_{k=0}^\infty a_k\, t^{-q-\frac12+k} \qquad\text{as}\quad t\to
0\notag
\end{align}
where the constant $c$ vanishes unless $q=0$. Proceeding as before
and using \eqref{e:asymp4},
\begin{align}\label{e:asymp5}
\int^\infty_{-\infty} \psi(1+2i\tau)h_{t,r}(\tau)\, d\tau\, \sim\,
\sum_{k=0}^\infty {a_k(r)}\, t^{-1+\frac{k}{2}}\,
+{b_k(r)}\,t^{-1+k}\, \log t \quad\text{as}\quad t\to 0
\end{align}
{where $a_k(r), b_k(r)$ depend only on $r$ and is of
$\mathrm{O}(r^2)$ for small $r>0$.}

For the fourth part of $\Tr_\p\big(\Dd_r e^{-t\Dd_r^2}\big)$, it
is also easy to get the following asymptotic expansion
\begin{equation}\label{e:non}
\int^\infty_{-\infty} h_{t,r}(\tau) \, d\tau \, \sim\,
\sum_{k=0}^\infty {a_k(r)}\, t^{-1+k} \qquad\text{as}\quad t\to 0
\end{equation}
where $a_k(r)$ depends only on $r$ and is of $\mathrm{O}(r^2)$ for
small $r>0$.

Now it is easy to see that the next term $\frac{\kappa_\t}{2}
h_{t,r}(0)$ contributes to the small time asymptotics by
\eqref{e:asymp3} {with the first term $a_1(r)t^{-\frac12}$}.

 By \eqref{e:J}, the integrand of the last integral of the
geometric side can be expressed by \begin{align}\label{e:digamma1}
\psi\big(\frac12+i\tau\big)+\psi\big(i\tau\big)-\psi\big(\frac{1+m}{2}+i\tau\big)-\psi\big(\frac{1-m}{2}+i\tau\big).
\end{align}
Using the following formulas about $\psi(z)$,
\begin{equation}\label{e:digamma2}
\psi(z+1)=\frac{1}{z}+\psi(z), \qquad
\psi(z)+\psi(z+\frac12)=2\big(\psi(2z)-\log 2\big), \end{equation}
the terms in \eqref{e:digamma1} can be rewritten as \begin{align*}
2\big(\psi(1+i\tau)&-\psi(1+2i\tau)\big)-\frac{1}{i\tau}+2\log 2\\
&- 4\Big(
\frac{1}{1+4\tau^2}+\frac{3}{3^2+4\tau^2}+\ldots+\frac{m-1}{(m-1)^2+4\tau^2}\Big).
\end{align*}
The terms in the first line gives us the asymptotics as
\eqref{e:asymp5}. The terms in the second line also can be handled
as in a similar way and we can show that these term gives us the
asymptotics
\begin{equation}\label{e:asymp6}
\sum_{k=0}^\infty {a_k(r)}\, t^{-\frac{1+k}{2}}
\qquad\text{as}\quad t\to 0
\end{equation}
{where $a_k(r)$ depends only on $r$ and is of $\mathrm{O}(r^2)$
for small $r>0$.} Combining \eqref{e:asymp3}, \eqref{e:asymp5},
\eqref{e:asymp6} and facts derived in the above, we obtain

\begin{theorem}\label{t:asymp} The small time asymptotics is given by
\begin{equation}\label{e:asymp-final}
\Tr_\p\big(\Dd_r e^{-t\Dd_r^2}\big) \, \sim\, \sum_{k=0}^\infty
{a_k(r)}\, t^{\frac{-3+k}{2}}\, +{b_k(r)}\,t^{-1+k}\, \log t
\quad\text{as}\quad t\to 0
\end{equation}
{where $a_k(r), b_k(r)$ depend only on $r$ and is of
$\mathrm{O}(r^2)$ for small $r>0$.} In particular, if
$\kappa^t=0$, it has the following simple form,
\begin{equation*}
\Tr_\p\big(\Dd_r e^{-t\Dd_r^2}\big) \, \sim\, \sum_{k=0}^\infty
{a_k(r)}\, t^{-\frac{3+k}{2}} \qquad\text{as}\quad t\to 0.
\end{equation*}
\end{theorem}

This theorem also immediately implies

\begin{theorem}\label{t:mer-eta-p} For a sufficiently small $r>0$, the function
$\eta_\p(\Dd_r,s)$ defined for $\Re(s)\gg 0$ has the meromorphic
extension over $\mathbb{C}$ and may have a double pole at $s=1$
and simple poles at $ -\mathbb{N}\cup \{0,2\}$. In particular, if
$\kappa^\t=0$, $\eta_\p(\Dd_r,s)$ may have only the simple poles
at $-\mathbb{N}\cup\{0,1,2\}$.
\end{theorem}

In the view of Theorem \ref{t:mer-eta-p}, it is natural to define
the principal part of the eta invariant of $\Dd_r$ by
\[
\eta_\p(\Dd_r):= \Big( \eta_\p(\Dd_r,s) -
\frac{r_0}{s}\Big)\Big|_{s=0}
\]
where $r_0$ is the residue of the simple pole of
$\eta_\p(\Dd_r,s)$ at $s=0$. Now let us consider the adiabatic
limit of $\eta_\p(\Dd_r)$ as $r\to 0$. For this, we need

\begin{proposition}\label{p:vanish} As $r\to 0$, $\Tr_\p\big(\Dd_r
e^{-t\Dd_r^2}\big)$ converges to $0$ for $t\in (0,\infty)$, and
$t^{\frac32}\Tr_\p\big(\Dd_r e^{-t\Dd_r^2}\big)$ converges to $0$
uniformly for $t\in [0,1]$.
\end{proposition}

\begin{proof}
From the expression of $h_{t,r}(\tau)$ in \eqref{e:htr}, we can
see
\[
|h_{t,r}(\tau)|\, \leq \, C\,\exp\big(-\frac{r^2}{4}t-4\tau^2
t\big) \sum_{m\in2\mathbb{Z}}e^{-c(m-1)^2(1+r^{-2})t}
\qquad\text{for small $r>0$}
\]
where $C,c$ are the positive constants that do not depend on
$r,\tau$. Hence, the integral \[\int^\infty_{-\infty} {\tau\tanh
(\pi \tau)}\, h_{t,r}(\tau) \, d\tau\] vanishes as $r\to 0$ by the
dominated convergence theorem. The same argument holds for other
terms defining $\Tr_\p\big(\Dd_r e^{-t\Dd_r^2}\big)$. Hence
$\Tr_\p\big(\Dd_r e^{-t\Dd_r^2}\big)$ converges  to $0$ as
$r\to0$. The uniform convergence of $t^{\frac32}\Tr_\p\big(\Dd_r
e^{-t\Dd_r^2}\big)$ follows from the following estimate
\begin{equation}\label{e:up-est}
\left|t^{\frac32}\Tr_\p\big(\Dd_r e^{-t\Dd_r^2}\big)\right| \leq C
r^2 \qquad \text{for} \ \ t\in [0,1],
\end{equation}
which also follows easily from \eqref{e:htr-est} and
\eqref{e:poisson}.
\end{proof}

Now we have

\begin{theorem}
\[
\lim_{r\to 0} \eta_\p(\Dd_r)=0.
\]
\end{theorem}

\begin{proof}
 Let us consider
\begin{align*}
\eta_\p(\Dd_r)= &\frac{1}{\sqrt{\pi}}\int^\infty_1 \,
t^{-\frac{1}{2}} \Tr_\p\big(\Dd_r e^{-t\Dd_r^2}\big)\, dt\\
&+ \Big(\frac{1}{\Gamma(\frac{s+1}{2})} \int^1_0 \,
t^{\frac{s-1}{2}} \Tr_\p\big(\Dd_r e^{-t\Dd_r^2}\big)\, dt -
\frac{r_0}{s}\Big)\Big|_{s=0}.
\end{align*}
For the integration over $[1,\infty)$, recalling that
$\Tr_\p\big(\Dd_r e^{-t\Dd_r^2}\big)$ is exponentially decaying as
$t\to\infty$, it is easy to see that this part vanishes as $r\to
0$ by Proposition \ref{p:vanish} and the dominated convergence
theorem. {By \eqref{e:asymp-final} the meromorphic extension of
the integral $\int^1_0\, \cdot\, dt$ has the following form for
$\Re(s)\geq -\epsilon$ with small $\epsilon>0$,
\begin{multline}\label{e:small-t}
 \int^1_0 \,
t^{\frac{s-1}{2}} \Tr_\p\big(\Dd_r e^{-t\Dd_r^2}\big)\, dt\, = \,
\frac{2a_0}{s-2} + \frac{2a_1}{s-1} - \frac{4b_0}{(s-1)^2}
+\frac{2a_2}{s}\\ + \frac{2a_3}{s+1}- \frac{4b_1}{(s+1)^2}
 + \int^1_0 \, t^{\frac{s-1}{2}} \Tr_\p^*\big(\Dd_r
e^{-t\Dd_r^2}\big)\, dt,
\end{multline}
where
\[
\Tr_\p^*\big(\Dd_r e^{-t\Dd_r^2}\big):= \Tr_\p\big(\Dd_r
e^{-t\Dd_r^2}\big) -a_0\, t^{-\frac32} - a_1\, t^{-1} - b_0
t^{-1}\log t - a_2\, t^{-\frac12} -a_3 - b_1\log t.
\]}
{By Theorem \ref{t:asymp},  all the coefficients $a_0,a_1,a_2,
a_3, b_0, b_1$ (as function of variable $r$) vanish as $r\to 0$.
Hence putting $s=0$ except the term $\frac{2a_2}{s}$, we can see
that $-a_0-2a_1-4b_1+2a_3-4b_1$ vanishes as $r\to 0$. For the last
integral with $s=0$ also vanishes as $r\to 0$ since
\[
\left|t^{-\frac12}\Tr_\p^*\big(\Dd_r e^{-t\Dd_r^2}\big)\right|
\leq C r^2 \qquad \text{for} \ \ t\in [0,1], \] which follows from
\eqref{e:up-est} and the coefficients $a_0,a_1,a_2, a_3, b_0, b_1$
vanish as order of $r^2$.} This completes the proof.

\end{proof}

%%%%%%%%%%%%%%%%%%%%%%%%%%%%%%%%%%%%%%%%%%%%%%%%%%%%%%%%%%%%%%%%%%
\section{Eta function of $\Dd_r$: Discrete series part}
\label{sec-etafnDiscreteseriespart}
%%%%%%%%%%%%%%%%%%%%%%%%%%%%%%%%%%%%%%%%%%%%%%%%%%%%%%%%%%%%%%%%%%%

In this section we study the discrete part of the eta function
$\eta_\d(\Dd_r,s)$ when $r>0$ is sufficiently small.

First, from Lemma \ref{l:dis-htr}, let us recall that $h_{t,r}(n)$
is given by $\lambda(n)$'s and $\lambda_\pm(n,m)$'s and we
decompose $\Tr_\d\big(\Dd_r e^{-t\Dd_r^2}\big)$ into the
corresponding two parts. Then we also have
\[
\eta_\d(\Dd_r,s)=\eta_\d^1(\Dd_r,s)+\eta_\d^2(\Dd_r,s)\qquad\text{for}
\quad \Re(s)\gg 0
\]
where
\[
\eta_\d^1(\Dd_r,s)= r^{s}\left(2g-2+{\kappa}\right)
\Big(-\sum_{k=1}^\infty
\frac{2(2k-1)}{(2k-1+\frac{r^2}{2})^s}\Big)
+r^s\kappa^\t\Big(\sum_{k=1}^\infty
\frac{2}{(2k-1+\frac{r^2}{2})^s}\Big)
\]
\begin{multline*}
\eta_d^2(\Dd_r,s)=\left(2g-2+{\kappa}\right)\Big(2\sum_{k=1}^\infty
(2k-1)\sum_{\ell\in k+\mathbb{N}}\lambda_+(2k,2\ell)^{-s}
-\lambda_-(2k,2\ell)^{-s}\Big)\\
-\kappa^\t\Big(2\sum_{k=1}^\infty \sum_{\ell\in
k+\mathbb{N}}\lambda_+(2k,2\ell)^{-s}
-\lambda_-(2k,2\ell)^{-s}\Big).
\end{multline*}
Here we used the fact
\[
{\mathrm{vol}(\Gamma\backslash G)} =
2\pi\left(2g-2+{\kappa}\right)
\]
where the volume of $\Gamma\backslash G$ is given w.r.t. the Haar
measure in \eqref{e:Haar} {(recall that the volume of the circle
$K/Z$ is normalized to be $1$)}.

 Now we investigate $\eta^1_\d(\Dd_r,s)$. Let us recall the
Hurwitz zeta function
\[
\zeta(s,a)=\sum_{k=0}^\infty (k+a)^{-s}
\]
which has a meromorphic extension to the whole $\mathbb{C}$ with a
simple pole at $s=1$. If we set
\[
\zeta_0(s,a)=\sum_{k=1}^\infty(2k-1+a)^{-s},
\]
then
\[
\zeta_0(s,a)=\zeta(s,a)-2^{-s}\zeta(s,\frac a2).
\]
By these definitions, for $\Re(s)\gg 0$,
\[
\eta^1_\d(\Dd_r,s)=2\left(2-2g-{\kappa}\right) r^s
\left(\zeta_0(s-1,\frac{r^2}{2})-\frac{r^2}{2}\zeta_0(s,\frac{r^2}{2})\right)
+2\kappa^\t r^s \zeta_0(s,\frac{r^2}{2}).
\]
The right hand side gives the meromorphic extension of
$\eta^1_\d(\Dd_r,s)$ over $\mathbb{C}$ with the simple poles at
$s=1,2$. We can also see that $\eta^1_\d(\Dd_r,s)$ is regular at
$s=0$ from this equality. Recalling
\[
\zeta(0,a)=\frac12- a, \qquad
\zeta(-1,a)=-\frac12\left(a^2-a+\frac16\right),
\]
we can see that
\[
\zeta_0(0,a)=-\frac {a}{2}, \qquad
\zeta_0(-1,a)=-\frac14\left(a^2-\frac13\right).
\]
Using these, we obtain
\[
\eta^1_\d(\Dd_r,0)=\left(2-2g-{\kappa}\right)\left(\frac16+\frac{r^4}{8}\right)-\kappa^\t\frac{r^2}{2}.
\]
Summarizing all these for $\eta^1_\d(\Dd_r,s)$, we have

\begin{proposition}\label{p:eta1} For a sufficiently small $r>0$,
the function $\eta_\d^1(\Dd_r,s)$ define for $\Re(s)\gg 0$ has the
meromorphic extension over $\mathbb{C}$ and has the simple poles
at $s=1,2$. The following equality holds,
\begin{equation}\label{e:eta1}
\lim_{r\to 0} \eta^1_\d(\Dd_r, 0)=\frac16
\left(2-2g-{\kappa}\right).
\end{equation}
\end{proposition}

To get the meromorphic extension of $\eta^2_\d(\Dd_r,s)$ over
$\mathbb{C}$, we rewrite this as follows,
\begin{equation*}
\eta^2_\d(\Dd_r,s)=2 \left(2g-2+{\kappa}\right) r^s f_r(s)
-2\kappa^\t r^s g_r(s).
\end{equation*}
Here
\begin{align*}
f_r(s)=&\sum_{\ell>k\geq1} (2k-1)\Big(
\big(q_r(k,\ell)-\frac{r^2}{2}\big)^{-s} -
\big(q_r(k,\ell)+\frac{r^2}{2}\big)^{-s}\Big),\\
g_r(s)=&\sum_{\ell>k\geq1}
\big(q_r(k,\ell)-\frac{r^2}{2}\big)^{-s} -
\big(q_r(k,\ell)+\frac{r^2}{2}\big)^{-s}
\end{align*}
where
\[
q_r(k,\ell)= \big((2\ell-1)^2(1+r^2)-r^2(2k-1)^2\big)^{\frac12}.
\]
Now we put $h_r(s)=\sum_{\ell>k\geq 1}(2k-1) q_r(k,\ell)^{-s}$
which can be written as
\begin{align*}
h_r(s) =\sum_{k\geq 1} (2k-1)^{1-s} \sum_{\ell>k}(2\ell-1)^{-s}
\Big(\frac{1+r^2}{(2k-1)^2}-\frac{r^2}{(2\ell-1)^2}\Big)^{-\frac
s2}.
\end{align*}
From this and the above analysis of $\zeta_0(s,0)$, we can see
that $h_r(s)$ is holomorphic for $\Re(s)>2$. For the meromorphic
extension of $h_r(s)$ over $\mathbb{C}$, we use the identity
$a^s=\exp(s\log(1+(a-1)))$ to get
\begin{align*}
&\Big(\frac{1+r^2}{(2k-1)^2}-\frac{r^2}{(2\ell-1)^2}\Big)^{-\frac
s2}\\ =& 1 -\frac s2
\big(\frac{1+r^2}{(2k-1)^2}-\frac{r^2}{(2\ell-1)^2}-1\big)+ \frac
s4 \big(\frac{1+r^2}{(2k-1)^2}-\frac{r^2}{(2\ell-1)^2}-1\big)^2
+\ldots .
\end{align*}
From this, we can see that $h_r(s)$ has the meromorphic extension
over $\mathbb{C}$ and may have the simple poles at
$s=2,1,0,-1,\ldots$ with the residues which are continuous w.r.t.
$r$. Using the following equality
\begin{align*}
f_r(s)= \left( r^2 s h_r(s+1) + r^6 \frac{s(s+1)(s+2)}{24}
h_r(s+3) + r^{10} \theta(s,r)\right)
\end{align*}
where $\theta(s,r)$ is regular at $s=0$ and is continuous at
$r=0$, we can conclude that $f_r(s)$ is regular at $s=0$ and the
limit of $f_r(0)$ as $r\to 0$ is trivial. In a similar way, we can
see that the same conclusion is true for $g_r(s)$. By all these
facts, we have

\begin{proposition}\label{p:eta2} For a sufficiently small $r>0$,
the function $\eta_\d^2(\Dd_r,s)$ defined for $\Re(s)\gg 0$ has
the meromorphic extension over $\mathbb{C}$ and may have the
simple poles at $s= -\mathbb{N}\cup \{1\}$. The following equality
holds,
\begin{equation}\label{e:eta2}
\lim_{r\to 0} \eta^2_\d(\Dd_r, 0)=0.
\end{equation}
\end{proposition}

By Proposition \ref{p:eta1}, \ref{p:eta2}, we can define
\[
\eta_\d(\Dd_r):=\eta_\d(\Dd_r,s)\big|_{s=0}=\eta^1_\d(\Dd_r,0)+\eta^2_\d(\Dd_r,0)
\]
and

\begin{theorem}For a sufficiently small $r>0$, the discrete part of the eta function
$\eta_\d(\Dd_r,s)$ has the meromorphic extension over $\mathbb{C}$
and may have the simple poles at $-\mathbb{N}\cup \{1,2\}$. The
following equality holds,
\begin{equation*}
\lim_{r\to 0}\eta_\d(\Dd_r)=\frac16 \left(2-2g-{\kappa}\right).
\end{equation*}
\end{theorem}

% ------------------------------------------------------------------------
\bibliographystyle{amsplain}
\def\cprime{$'$} \def\polhk#1{\setbox0=\hbox{#1}{\ooalign{\hidewidth
  \lower1.5ex\hbox{`}\hidewidth\crcr\unhbox0}}}
\providecommand{\bysame}{\leavevmode\hbox
to3em{\hrulefill}\thinspace}
\providecommand{\MR}{\relax\ifhmode\unskip\space\fi MR }
\providecommand{\MRhref}[2]{%
  \href{http://www.ams.org/mathscinet-getitem?mr=#1}{#2}
} \providecommand{\href}[2]{#2}

\end{document}